\documentclass[12pt]{iopart}

\usepackage{iopams}  
\usepackage{algorithm}
\usepackage{algpseudocode}
\usepackage{amssymb,graphicx,graphics}
\usepackage{graphicx}
\usepackage{amsthm, soul} 
\usepackage{iopams}
\expandafter\let\csname equation*\endcsname\relax
\expandafter\let\csname endequation*\endcsname\relax
\usepackage{amsmath}

\usepackage{enumitem}
\usepackage{comment}
\usepackage{booktabs}

\usepackage{verbatim}
\usepackage{color}
\usepackage[colorinlistoftodos]{todonotes}
\usepackage{tikz}
\usepackage{nicefrac}

\newtheorem{proposition}{Proposition}

\theoremstyle{definition} 
\newtheorem{remark}{Remark}

\definecolor{dgreen}{rgb}{0,0.20,0}

\begin{document}

\title[Solving Implicit Inverse Problems with Homotopy-Based Regularization Path]{Solving Implicit Inverse Problems with Homotopy-Based Regularization Path}%Regularizing IIP via the Homotopy Method

\author{Davide Parodi$^{1,2}$, Federico Benvenuto$^{2}$, Sara Garbarino$^{2,3}$, Michele Piana$^{2,3}$}
\address{$^1$ Universit\`a Campus Bio-Medico di Roma}
\address{$^2$ MIDA, Dipartimento di Matematica, Universit\`a di Genova, via Dodecaneso 35 16146 Genova, Italy}
\address{$^3$ Life Science Computational Laboratory, IRCCS Ospedale Policlinico San Martino, Genova, Italy}

\date{May 2024}

\begin{abstract}
Implicit inverse problems, in which noisy observations of a physical quantity are used to infer a nonlinear functional applied to an associated function, are inherently ill-posed and often exhibit non-uniqueness of solutions. Such problems arise in a range of domains, including the identification of systems governed by Ordinary and Partial Differential Equations (ODEs/PDEs), optimal control, and data assimilation. Their solution is complicated by the nonlinear nature of the underlying constraints and the instability introduced by noise.
In this paper, we propose a homotopy-based optimization method for solving such problems. Beginning with a regularized constrained formulation that includes a sparsity-promoting regularization term, we employ a gradient-based algorithm in which gradients with respect to the model parameters are efficiently computed using the adjoint-state method. Nonlinear constraints are handled through a Newton–Raphson procedure.
By solving a sequence of problems with decreasing regularization, we trace a solution path that improves stability and enables the exploration of multiple candidate solutions. The method is applied to the latent dynamics discovery problem in simulation, highlighting performance as a function of ground truth sparsity and semi-convergence behavior.

\end{abstract}


\section{Introduction}
Implicit inverse problems \cite{mosegaard1999probabilistic,nino2018local} require the determination of a non-linear functional applied to a function representing a physical quantity, when only noisy measurements of this quantity are available. 
When the functional is modeled by a parametric form, the implicit inverse problem becomes the one to optimize its parameters. 
Due to the inherent complexity of the problem and the nonlinear nature of their constraints, implicit inverse problems are typically strongly ill-posed and, in particular,  are characterized by non-uniqueness of the solution. 

Typical examples of these problems are found in Ordinary or Partial Differential Equation constrained optimization, such as equation learning, shape optimization, optimal control and data assimilation problems.

The aim is to determine the parameters of an ODE or PDE system such that the corresponding solution best fits the observed data \cite{brunton2016discovering, chen2021physics, heinonen2018learning, lorenzi2018constraining, dondelinger2013ode, williams2006gaussian, rielly2025mock}. 
This involves the optimization of model parameters to match the observation and is often applied in latent dynamics discovery, where the underlying dynamics governing a system need to inferred from data.
These problems are common in fields like geophysics and biological modeling, where the governing equations of a system are either unknown or poorly understood, and must be learned from data \cite{plessix2006review}. 

In optimal control, the goal is to determine a control function that minimizes an objective function while satisfying constraints imposed by a dynamic system of ODEs or PDEs. This type of problem occurs in a wide range of applications, from robotics (where the goal is to optimize control inputs to achieve specific movements) to medicine (for example, optimizing drug dosages) and finance (where the objective might be to minimize financial risk) \cite{bergounioux2019position, alexanderian2021optimal}.

Data assimilation involves integrating observational data with mathematical models to obtain the best possible representation of a physical system.  This is particularly important in areas such as weather forecasting and climate modeling, where models like the Navier-Stokes equations are combined with real-world data (e.g., satellite measurements) to predict future states of the system. In methods such as variational data assimilation, the goal is to optimize the model parameters to align with the observed data, often in high-dimensional spaces \cite{baayen2019overview, keller2024ai, burger2020data}.

% REF DI PAPER DA AGGIUNGERE SOPRA, NEI CONTESTI SPECIFICI ODE, PDE, ECC ECC) \cite{dondelinger2013ode, brunton2016discovering, chen2021physics, heinonen2018learning, lorenzi2018constraining, williams2006gaussian, rielly2025mock, burger2020data}.

Given the broad range of problems that can be considered as implicit inverse problems, the proposed solution methods are typically problem-specific, and the methods are often not framed within the broader context of implicit inverse problems, with few exceptions \cite{rothermel2021solving, migorski2019inverse, bergounioux2019position, alexanderian2021optimal}. Indeed, the strong ill-posedness of these problems often forces the development of solution methods that are restricted to the specific application.

In this paper, we propose a homotopy method \cite{watson1989modern, baayen2019overview} to follow the regularization path that can potentially be applied to all specific problems that can be framed as implicit inverse problems.
The method exploits a variational approach by minimizing a functional that combines a data-fitting term and a penalty term to enforce the stability of the solution.
For any fixed value of the regularization parameter, we solve the problem using a gradient algorithm that leverages the adjoint state method \cite{plessix2006review} for computing large Jacobians and the Newton-Raphson technique \cite{Dedieu2015} for solving nonlinear forward problems at each iteration. These consolidated numerical tools ensure the robustness of the algorithm.
Our strategy is to start with large values of the regularization parameter, ensuring that the solution exists and is unique.
Then, we iteratively solve the problem with decreasing values of the regularization parameter to obtain a complete regularization path \cite{friedman2010regularization}.
Moreover, we apply the proposed method to a latent dynamics discovery problem, namely an ODE learning problem, where the forward model is given by a an autonomous differential equation characterized by few important terms that govern the dynamics. 
We show how this method can recover the important terms of the dynamics in the case the differential equations contains only periodic functions. 
As in the case of classical linear inverse problems, we show the well known semi-convergent behavior of the regularized solution family.

The plan of the paper is as follows. Section 2 provides the mathematical framework of the numerical method, which is formulated in Section 3. Section 4 describes an application related to the problem of learning ODEs. Section 5 concludes the paper with a summary of our finding and conclusions.

\section{Mathematical framework}
Let $\mathcal{H}$ be a complete Hilbert space and $u \in \mathcal{H}^U$ such that
$$
t \rightarrow u(t) = [u_0(t), u_1(t),...,u_{U-1}(t)] \in \mathbb{R}^U
$$
where $u_h \in \mathcal{H}$ for each $h=0,...,U-1$. Let
$$
    \mathcal{F}: \mathcal{H}^U \times \mathbb{R}^M \rightarrow \mathcal{H}^U
$$
$$
(u, \textbf{m}) \rightarrow \mathcal{F}(u, \textbf{m})
$$
be an implicit parametric functional, where the parametric form is known but the values of the entries of $\textbf{m} \in \mathbb{R}^M$ are unknown. The relationship between $u$ and $\textbf{m}$ can be expressed as
\begin{equation}\label{State equation}
\mathcal{F}(u, \textbf{m}) = 0. 
\end{equation} 

We now introduce two vectors, $\textbf{u}, \textbf{d} \in \mathbb{R}^{T \times U}$, such that

\begin{equation}\label{vector u}
u_{k,h} = u_h(t_k) \ \ \ \ k=0,...,T-1 \ \ ; \ \  h=0,...,U-1,
\end{equation}
and 
\begin{equation}\label{vector d}
    d_{k,h} = u_{k,h} + \eta,
\end{equation}
where $\eta$ is sampled from a Gaussian distribution $\mathcal{N}(\textbf{0}, \sigma)$.
Thus, the parametric implicit inverse problem addressed in this study is to determine $\textbf{m}$ in \eqref{State equation}, given a discretized, noisy version $\textbf{d}$ of $u$ from \eqref{vector d} and the known parametric form of $\mathcal{F}$.

To address the intrinsic ill-posedness of the problem \cite{hadamard1902problemes}, we introduce a regularized objective function $k_{\textbf{d}, \alpha}: \mathbb{R}^{T \times U} \times \mathbb{R}^{M} \rightarrow \mathbb{R}$ defined as
\begin{equation}\label{Objective Functional}
(\textbf{u}, \textbf{m}) \rightarrow k_{\textbf{d}, \alpha}(\textbf{u}, \textbf{m}) = h_\textbf{d}(\textbf{u}) + \alpha g(\textbf{m}) \ \ ,
\end{equation}
where the choice 
\begin{equation}\label{data loss}
h_{\textbf{d}}(\textbf{u}) = \| \textbf{u} - \textbf{d} \|_2^2
\end{equation}
is consistent with the Gaussian noise affecting the measurements \eqref{vector d}. In \eqref{Objective Functional}, $\alpha > 0$ is the regularization parameter, and $g: \mathbb{R}^M \rightarrow \mathbb{R}$ is the regularization function \cite{tikhonov1977solutions, rudin1992nonlinear} [REF!!!!].

We now define a discretized version of the implicit parametric functional in \eqref{State equation}:
\begin{equation}\label{discretized Implicit Functional}
    (\textbf{u}, \textbf{m}) \rightarrow \textbf{F}(\textbf{u}, \textbf{m})
\end{equation}
so that equation \eqref{State equation} becomes
\begin{equation}\label{State equation discrete}
    \textbf{F}(\textbf{u}, \textbf{m}) = \textbf{0}.
\end{equation}
We denote by $\textbf{u}_{\textbf{m}} \in \mathbb{R}^{T \times U}$ the solution to \eqref{State equation discrete} for a given set of parameters $\textbf{m}$.

From an operational standpoint, the solution of the discretized, regularized parametric implicit inverse problem is given by
\begin{equation}\label{Minimum regluarized constrained problem discrete}
        \textbf{m}_{\alpha} \in \arg \min_{\textbf{m} \in \mathbb{R}^M} k_{\textbf{d}, \alpha}(\textbf{u}_{\textbf{m}}, \textbf{m})
\end{equation}
where $\textbf{d}$ is the input data vector of the problem.
The solution to the corresponding non-regularized discretized inverse problem is
\begin{equation}\label{Minimum constrained problem discrete}
    \textbf{m}^* \in \arg \min_{\textbf{m} \in \mathbb{R}^M} h_\textbf{d}(\textbf{u}_{\textbf{m}}),
\end{equation}
which is obtained by setting $\alpha = 0$ in \eqref{Minimum regluarized constrained problem discrete}.

\section{Homotopy-based regularization} 
We propose a homotopy-based optimization scheme to regularize the implicitly constrained inverse problem defined by equation \eqref{Minimum constrained problem discrete}. The method is structured in two conceptual layers:
\begin{itemize}
    \item Inner loop – adjoint-based optimization: for a fixed regularization parameter $\alpha > 0$, we solve the regularized problem via a gradient-based algorithm using the adjoint state method. This yields a sequence of approximations to the parameter vector $\textbf{m}_\alpha$.
    \item Outer loop – homotopy continuation: We trace the regularization path by progressively decreasing $\alpha$ and using the solution obtained for $\alpha_l$ to initialize the optimization for $\alpha_{l+1}$ (a warm restart). This allows us to gradually approach the solution $\textbf{m}^*$ of the non-regularized problem while maintaining numerical stability.
\end{itemize}

\subsection{Inner loop: adjoint-based optimization}
\label{subsec:adjoint_method}

For a given $\alpha > 0$, we solve the regularized inverse problem by minimizing the Lagrangian:
\begin{equation}\label{Lagrangian}
\mathcal{L}_{\textbf{d}, \alpha}(\textbf{u}, \textbf{m}, \boldsymbol{\lambda}) = h_\textbf{d}(\textbf{u}) + \alpha g(\textbf{m}) - \langle \boldsymbol{\lambda}, \textbf{F}(\textbf{u}, \textbf{m}) \rangle_{T \times U},
\end{equation}
where $\textbf{F}(\textbf{u}, \textbf{m}) = 0$ encodes the implicit constraint.

We recall the following result characterizing stationarity of the regularized problem:
\proposition{Necessary condition for the solution of the discretized regularized implicit inverse problem.}\label{prop:necessary condition}
Given the discretized regularized parametric implicitly constrained minimum problem \eqref{Minimum regluarized constrained problem discrete} associated to functional $\textbf{F}$, one of its solution $\textbf{m}_{\alpha}$ and the associated vector $\textbf{u}_{\textbf{m}_{\alpha}}$, then there exists $\boldsymbol{\lambda}_{\textbf{m}_{\alpha}} \in \mathbb{R}^{T \times U}$ such that all the following conditions are fulfilled:
\begin{enumerate}
    \item $\boldsymbol{\lambda}_{\textbf{m}_{\alpha}} = \left ( \left [ \dfrac{\partial \textbf{F}}{\partial \textbf{u}} \right]^T \right )^{\dag} \dfrac{\partial h_{\textbf{d}}}{\partial \textbf{u}} \Big|_{(\textbf{u}_{\textbf{m}_{\alpha}}, \textbf{m}_{\alpha})}$;
    \item $\alpha \dfrac{\partial g}{\partial \textbf{m}} - \langle \boldsymbol{\lambda}, \dfrac{\partial \textbf{F}}{\partial \textbf{m}} \rangle_{T \times U} \Big|_{(\textbf{u}_{\textbf{m}_{\alpha}}, \textbf{m}_{\alpha}, \boldsymbol{\lambda}_{\textbf{m}_{\alpha}})} = \textbf{0}$. 
\end{enumerate}

\proof{
Let us set $\textbf{u}_{\textbf{m}_{\alpha}} = \textbf{u}_{\alpha}$ and $\boldsymbol{\lambda}_{\textbf{m}_{\alpha}} = \boldsymbol{\lambda}_{\alpha}$. We can define the Lagrangian of the problem as
\begin{equation}
        \mathcal{L}_{\textbf{d}, \alpha}(\textbf{u}, \textbf{m}, \boldsymbol{\lambda}) = h_\textbf{d}(\textbf{u}) + \alpha g(\textbf{m}) - \langle \boldsymbol{\lambda}, \textbf{F}(\textbf{u}, \textbf{m}) \rangle_{T \times U},
\end{equation}
where $
\langle \boldsymbol{\lambda}, \textbf{F}(\textbf{u}, \textbf{m}) \rangle_{T \times U} = \sum_{t=0}^{T-1} \sum_{h=0}^{U-1} \lambda_{t,h} \textbf{F}(\textbf{u}, \textbf{m})_{t,h}$. According to Lagrange's multiplier theorem \cite{bertsekas2014constrained}, if $\textbf{m}_{\alpha}$ is a minimum for the discretized \eqref{Minimum regluarized constrained problem discrete}, then there must exist $\boldsymbol{\lambda}_{\alpha} \in \mathbb{R}^{T \times U}$ such that $(\textbf{u}_{\alpha}, \textbf{m}_{\alpha},\boldsymbol{\lambda}_{\alpha})$ is a stationary point for the Lagrangian:
\begin{equation}\label{Lagrangian stationary point}
    \nabla \mathcal{L}_{\textbf{d}, \alpha} (\textbf{u}_{\alpha}, \textbf{m}_{\alpha}, \boldsymbol{\lambda}_{\alpha}) = \textbf{0}
\end{equation}
This implies that all partial derivatives must be equal to 0 at $(\textbf{u}_{\alpha}, \textbf{m}_{\alpha}, \boldsymbol{\lambda}_{\alpha})$:
\begin{equation}\label{lagrangian p d lambda}
    \dfrac{\partial \mathcal{L}_{\textbf{d}, \alpha}}{\partial \boldsymbol{\lambda}}(\textbf{u}_{\alpha}, \textbf{m}_{\alpha}, \boldsymbol{\lambda}_{\alpha}) = \textbf{F}(\textbf{u}_{\alpha},\textbf{m}_{\alpha}) = \textbf{0};
\end{equation}
\begin{equation}\label{lagrangian p d u}
    \dfrac{\partial \mathcal{L}_{\textbf{d}, \alpha}}{\partial \textbf{u}}(\textbf{u}_{\alpha}, \textbf{m}_{\alpha}, \boldsymbol{\lambda}_{\alpha}) = \dfrac{\partial h_{\textbf{d}}}{\partial \textbf{u}} (\textbf{u}_{\alpha}) - \langle \boldsymbol{\lambda}_{\alpha}, \dfrac{\partial \textbf{F}}{\partial \textbf{u}} (\textbf{u}_{\alpha}, \textbf{m}_{\alpha}) \rangle_{T \times U} = \textbf{0};
\end{equation}
\begin{equation}\label{lagrangian p d m}
    \dfrac{\partial \mathcal{L}_{\textbf{d}, \alpha}}{\partial \textbf{m}}(\textbf{u}_{\alpha}, \textbf{m}_{\alpha}, \boldsymbol{\lambda}_{\alpha}) = \alpha \dfrac{\partial g}{\partial \textbf{m}}(\textbf{m}_{\alpha}) - \langle \boldsymbol{\lambda}_{\alpha}, \dfrac{\partial \textbf{F}}{\partial \textbf{m}}(\textbf{u}_{\alpha}, \textbf{m}_{\alpha})\rangle_{T \times U} = \textbf{0}.
\end{equation}
Equation \eqref{lagrangian p d lambda} is equivalent to equation \eqref{State equation discrete} with parameters $\textbf{m}_{\alpha}$. Equation \eqref{lagrangian p d u} implies that $\boldsymbol{\lambda}_{\alpha}$ must solve
\begin{equation}\label{adjoint State System}
           \left [ \dfrac{\partial \textbf{F}}{\partial \textbf{u}} (\textbf{u}_{\alpha}, \textbf{m}_{\alpha}) \right]^T \boldsymbol{\lambda}_{\alpha} = \dfrac{\partial h_{\textbf{d}}}{\partial \textbf{u}}(\textbf{u}_{\alpha}),
\end{equation}
which means that
\begin{equation}\label{adjoint State Variable}
    \boldsymbol{\lambda}_{\alpha} = \left ( \left [ \dfrac{\partial \textbf{F}}{\partial \textbf{u}} \right]^T \right )^{\dag} \dfrac{\partial h_{\textbf{d}}}{\partial \textbf{u}} \Big|_{(\textbf{u}_{\alpha}, \textbf{m}_{\alpha})}
\end{equation}
Replacing $\textbf{u}_{\alpha}$ and $\boldsymbol{\lambda}_{\alpha}$ obtained by solving \eqref{lagrangian p d lambda} and \eqref{adjoint State System} respectively in  
\eqref{lagrangian p d m}, we demonstrate point (ii) of the preposition.
}

\remark{The adjoint state method consists of computing the gradient of the Lagrangian \eqref{Lagrangian} with respect to the unknown parameters $\textbf{m}$ by computing $\boldsymbol{\lambda}_{\alpha}$, called adjoint state variables obtained by solving the linear system \eqref{adjoint State System}, called adjoint state system.}

\remark{ This method can be applied to any $h_\textbf{d}: \mathbb{R}^{T \times U} \rightarrow \mathbb{R}$ such that:
    \begin{itemize}
        \item $h_\textbf{d}(\textbf{u}) \ge 0$;
        \item $h_\textbf{d}(\textbf{u}) = 0 \iff u_h(t_k) = d_{k,h} $ for each $k=0,...,T-1, \; h=0,..., U-1$.
    \end{itemize}
Under these conditions $h_{\textbf{d}}$ can act as an empirical loss function \cite{rosasco2004loss}. 

Thanks to Proposition \ref{prop:necessary condition}, we can easily derive a gradient descent algorithm where each iteration consists of:
\begin{enumerate}
    \item find solution $\textbf{u}_{\textbf{m}_{\alpha}^{(i)}}$ of equation \eqref{State equation discrete};
    \item find $\boldsymbol{\lambda}_{\textbf{m}_{\alpha}^{(i)}}$ solution of \eqref{adjoint State System} computed in $\textbf{u}_{\textbf{m}_{\alpha}^{(i)}}, \textbf{m}_{\alpha}^{(i)}$;
    \item compute $\dfrac{\partial \mathcal{L}_{\textbf{d}, \alpha}}{\partial \textbf{m}} (\textbf{u}_{\textbf{m}_{\alpha}^{(i)}}, \textbf{m}_{\alpha}^{(i)}, \boldsymbol{\lambda}_{\textbf{m}_{\alpha}^{(i)}})$ as in \eqref{lagrangian p d m};
    \item update parameters $\textbf{m}_{\alpha}^{(i)}$ through gradient descent algorithm
    \begin{equation}\label{Gradient Step}
        \textbf{m}_{\alpha}^{(i+1)} = \textbf{m}_{\alpha}^{(i)} - \tau \dfrac{\partial \mathcal{L}_{\textbf{d}, \alpha}}{\partial \textbf{m}} (\textbf{u}_{\textbf{m}_{\alpha}^{(i)}}, \textbf{m}_{\alpha}^{(i)}, \boldsymbol{\lambda}_{\textbf{m}_{\alpha}^{(i)}})
    \end{equation}
    with $\tau >0$ the gradient step.
\end{enumerate}

\remark{As $\textbf{F}$ is non-linear, neither the forward nor adjoint equations are solved to full precision. Inexact solutions are used to reduce computational cost and mitigate instability by means of Newton-Raphson (NR) method \cite{Dedieu2015}.}

\remark{Step (ii) of the optimization algorithm is solved by using Landweber technique (LW) \cite{hanke1995convergence, scherzer1995convergence}.}

\remark{The optimization algorithm consists of an iterative update of $(\textbf{u}_{\textbf{m}_{\alpha}^{(i)}}, \textbf{m}_{\alpha}^{(i)}, \boldsymbol{\lambda}_{\textbf{m}_{\alpha}^{(i)}})$ for $i\ge0$. These quantities are never computed exactly but approximated in order to avoid numerical instability. 
At the end of each iteration of the inner loop, early stopping \cite{yao2007early} is applied. Internal optimization algorithms (NR and LW) are never brought to convergence either.}

The output of this iterative algorithm corresponding to regularization parameter $\alpha$ is the set of quantities $(\textbf{u}_{\alpha}, \textbf{m}_{\alpha}, \boldsymbol{\lambda}_{\alpha})$.

\subsection{Outer loop: homotopy continuation}

To recover the solution of the original (non-regularized) problem, we construct a regularization path:

\begin{equation}\label{regularization parameter}
    \boldsymbol{\alpha} = [\alpha_0, \alpha_1,..., \alpha_L]
\end{equation} such that $\alpha_0 > \alpha_1 >...> \alpha_L >0$.
For each $\alpha_l$, we solve the regularized problem using the adjoint-based algorithm above, initializing from the solution at $\alpha_{l-1}$:

\begin{enumerate}
    \item Compute $\textbf{m}_{\alpha_l}$ using the algorithm from Subsection~\ref{subsec:adjoint_method};
    
    \item Use $(\textbf{u}_{\alpha_l}, \textbf{m}_{\alpha_l}, \boldsymbol{\lambda}_{\alpha_l})$ to initialize the optimization for $\alpha_{l+1}$ \cite{kukavcka2017regularization}.
\end{enumerate}

\remark{This strategy enables stable convergence even in ill-posed settings. It ensures that the algorithm follows a continuous path of minimizers, which typically remain in the basin of attraction of a desired local solution.}

\subsection{Algorithmic Summary}

The complete procedure is summarized in Algorithm~\ref{alg:RASM}. Each outer loop iteration corresponds to a step in the homotopy, while each inner loop performs adjoint-based optimization for fixed $\alpha_l$.

\begin{algorithm}[ht]
\caption{Homotopy-based Adjoint Optimization}\label{alg:RASM}
\begin{algorithmic}[1]
\State Initialize $\tau > 0$, $\textbf{m}^{(0,0)}$
\For{$l = 0$ to $L$}
\State Set $\alpha \gets \alpha_l$
\For{$i = 1$ to $N_{\textit{max}}$}
\State Solve $\textbf{F}(\textbf{u}^{(l,i)}, \textbf{m}^{(l,i)}) = 0$ for $\textbf{u}^{(l,i)}$ (NR method)
\State Solve adjoint equation for $\boldsymbol{\lambda}^{(l,i)}$ (LW method)
\State Compute $\nabla_{\textbf{m}} \mathcal{L}_{\textbf{d}, \alpha}$
\State Update $\textbf{m}^{(l,i)}$ and apply early stopping
\EndFor
\State Set $\textbf{m}^{(l+1,0)} \gets \textbf{m}^{(l)} \gets \textbf{m}^{(l, N{\textit{iter}})}$
\EndFor
\State \Return The entire regularization path: $(\textbf{u}^{(l)}, \textbf{m}^{(l)}, \boldsymbol{\lambda}^{(l)})$ for $l=0,...,L$
\end{algorithmic}
\end{algorithm}

\section{Application to latent dynamic discovery problem}

In this section, we apply the proposed algorithm to the problem of discovering latent dynamics, where the goal is to reconstruct the unknown dynamics of a set of continuous physical quantities, denoted by $u$ (the state variables), from their temporal samples measured at times $t_0,\dots,t_{T-1}$ and collected in the data vector $\textbf{d}$. We assume that the state variables are periodic with respect to the temporal samples, meaning $\textbf{u}(t_0) = \textbf{u}(t_{T-1})$. Inspired by the formulation in \cite{heinonen2018learning}, the governing equation \eqref{State equation} simplifies to:
\begin{equation}\label{State equation Expl}
    \mathcal{F}(u, \textbf{m}) = \dot{u} - f(u, \textbf{m}) = 0, 
\end{equation}
where the function
$$  f: \mathcal{H}^U \times \mathbb{R}^M \rightarrow \mathcal{H}^U $$
$$ (u, \textbf{m}) \rightarrow f(u, \textbf{m})$$
describes the system's dynamics. 
Following \cite{brunton2016discovering, chen2018neural, heinonen2018learning}, we assume that the latent dynamics can be expressed as a linear combination of basis functions:
\begin{equation}\label{basis function}
    \Phi(u) = \{\phi_0(u), \phi_1(u),..., \phi_{D-1}(u)  \}
\end{equation}
where, for each $j=0,...,D-1$
$$\phi_j: \mathcal{H}^U \rightarrow \mathcal{H}^U$$
$$  u \rightarrow \phi_j(u) = [\phi_j(u)_0, \phi_j(u)_1,...,\phi_j(u)_{U-1}],$$
is a vector-valued function with components 
$\phi_j(u)_h \in \mathcal{C}^{1}(\mathcal{H}^U, \mathcal{H})$ for each $h=0,...,U-1$.
The dynamics are then modeled as:
\begin{equation}\label{Derivative function sparse}
    f_h(u, \textbf{m}) =  \sum_{j=0}^{D-1} \sum_{h'=0}^{U-1} m_{h'+Uj,h} \phi_j(u)_{h'}
\end{equation} 
for each component $h=0,...,U-1$, where $\textbf{m} \in \mathbb{R}^{(D \times U) \times U}$ is the unknown parameter matrix.
This formulation is consistent with the working hypothesis of knowing the form of implicit parametric functional $\mathcal{F}$. This formulation is consistent with the hypothesis that the parametric form of the operator $\mathcal{F}$ is known.
The corresponding dynamics can be described by the Cauchy problem:

\begin{equation}\label{Cauchy problem}
    \begin{cases}
    u(t_0) = u_0, \\
        \dot{u}_h -  \sum_{j=0}^{D-1} m_{jh} \phi_j(u)_h = 0 & \text{for each $h=0,...,U-1$}
    \end{cases}
\end{equation} 
given an initial condition $u(t_0) = u_0 \in \mathbb{R}^U$.

By evaluating $u$ and $\Phi(u)$ at the discrete time samples $t_0, \dots, t_{T-1}$, we obtain the matrix $\boldsymbol{\Phi} \in \mathbb{R}^{T \times (D \times U)}$, with components $\phi_j(\textbf{u})_h \in \mathbb{R}^T$. The discretized form of the functional $\textbf{F}$ in \eqref{State equation discrete}, under the assumption of \eqref{State equation Expl}, becomes:
\begin{equation}\label{State equation discrete latent dynamics} 
    \textbf{F}(\textbf{u}, \textbf{m}) = \dot{\textbf{u}} - \langle \boldsymbol{\Phi}^T(\textbf{u}), \textbf{m} \rangle = \textbf{0}
\end{equation}
where $\dot{\textbf{u}} \in \mathbb{R}^{T \times U}$ contains the discrete time derivatives of the state variables.

To promote sparsity in the learned dynamics and reduce model complexity, we apply $L_1$-regularization to the coefficient matrix $\textbf{m}$. In order to enable the use of the adjoint method, we use a smooth approximation:
\begin{equation}\label{regularization term}
    g(\textbf{m}) = \sum_{p=1}^M \sqrt{m_{p} + \epsilon^2},
\end{equation}
with $\epsilon > 0$ small, as suggested in \cite{defrise2011algorithm, vogel2002computational}. While this yields nearly sparse solutions, true sparsity is enforced by a hard thresholding step \cite{german2011nonlinear} applied after each parameter update:
\begin{equation}\label{Hard Thresh}
HT(m^{(l,i)}_p, \alpha_l) = 
    \begin{cases}
    0, & \text{if $| m^{(l,i)}_p | < \dfrac{\alpha_l}{2}$} \\
    m^{(l,i)}_p & \text{otherwise}
    \end{cases}
\end{equation} 
where $\alpha_l$ is the regularization parameter. This step is implemented at line 7 of Algorithm \ref{alg:RASM}.

\subsection{Results}
We evaluated the proposed method in a synthetic setting for the one-dimensional case ($U=1$) using basis function:
\begin{equation}\label{eq:basis}
\Phi(u) = (\cos(u), \cos(2u), \cos(3u), \cos(4u), \cos(5u), \cos(6u)),
\end{equation}
which includes trigonometric terms that can be used to approximate any periodic function.
We generated synthetic datasets using the ground truth coefficient vectors:
\begin{equation}\label{b1-m1}
\textbf{m}_1 = (1, -1, 0, 0, 0, 0)
\end{equation}
\begin{equation}\label{b1-m2}
\textbf{m}_2 = (-1.5, 1.5, -1.5, 1, -1, 0)
\end{equation}
and solved the forward Cauchy problem \eqref{Cauchy problem} with $u(t_0) = 0.2$ using the Python ODE solver \textit{solveivp} from the \textit{scipy.integrate} package. White Gaussian noise with varying standard deviations was added to simulate measurement noise.
Parameter settings are:
\begin{itemize}
    \item $\textbf{t}_k = \frac{2 \pi k}{T-1}$ for $k=0,...,T-1$ and $T=100$ the time samples;
    \item $\eta \sim N(0, \sigma)$ the additive white Gaussian noise, where we considered three values for $\sigma$, i.e., $\sigma=0.01, 0.1, 0.2$.
    \item $\boldsymbol{\alpha}$ a vector of $L=100$ logarithmically equi-spaced points in the interval $[ 10^{0}, 10^{-6}]$;
    \item $\tau = 10^{-3}$ gradient step;
    \item $N_{max} = 1000$ maximum number of Lagrangian optimization algorithm with the adjoint state method;
    \item $R_{max} = 50$ maximum number of NR iterations;
    \item $L_{max} = 100$ maximum number of LW iterations;
    \item $N_{ES} = 5$ number of checked loss function values for early stopping technique.  
\end{itemize}

For each ground truth $\textbf{m}$ and noise level $\sigma$, we generated $n=20$ realizations of the noise and applied the algorithm. For each trial $q = 1,...,n$, we define
\begin{equation}\label{eq:best-alpha}
    \alpha_q^* = \arg \min_{\alpha^{(l)}, l=0,...,99} \frac{\|{\textbf{m}_q^{(l)} - \textbf{m}\|_2}}{\| \textbf{m}\|_2}
\end{equation}
To keep the notation simple, we will omit the noise level dependency. This will be clear and explicit in every situation.

We define $\textbf{m}_q^*$ the set of parameters associated with the best regularization parameter $ \alpha_q^*$ \eqref{eq:best-alpha}. $\textbf{u}^*_{q}$ is the solution curve of the Cauchy problem \eqref{Cauchy problem} associated with $\textbf{m}_q^*$. The relative mean square errors and standard deviations are computed with respect to the $q=1,..., n=20$ trials.

Figure \ref{fig:synthetic-data-b1} shows the generated synthetic data for $\textbf{m}_1$ and $\textbf{m}_2$ under the three noise levels $\boldsymbol{\sigma} = [0.01, 0.1, 0.2]$.

\begin{figure}[H]
    \centering
    \begin{tabular}{c c c}
    \includegraphics[scale=0.29]{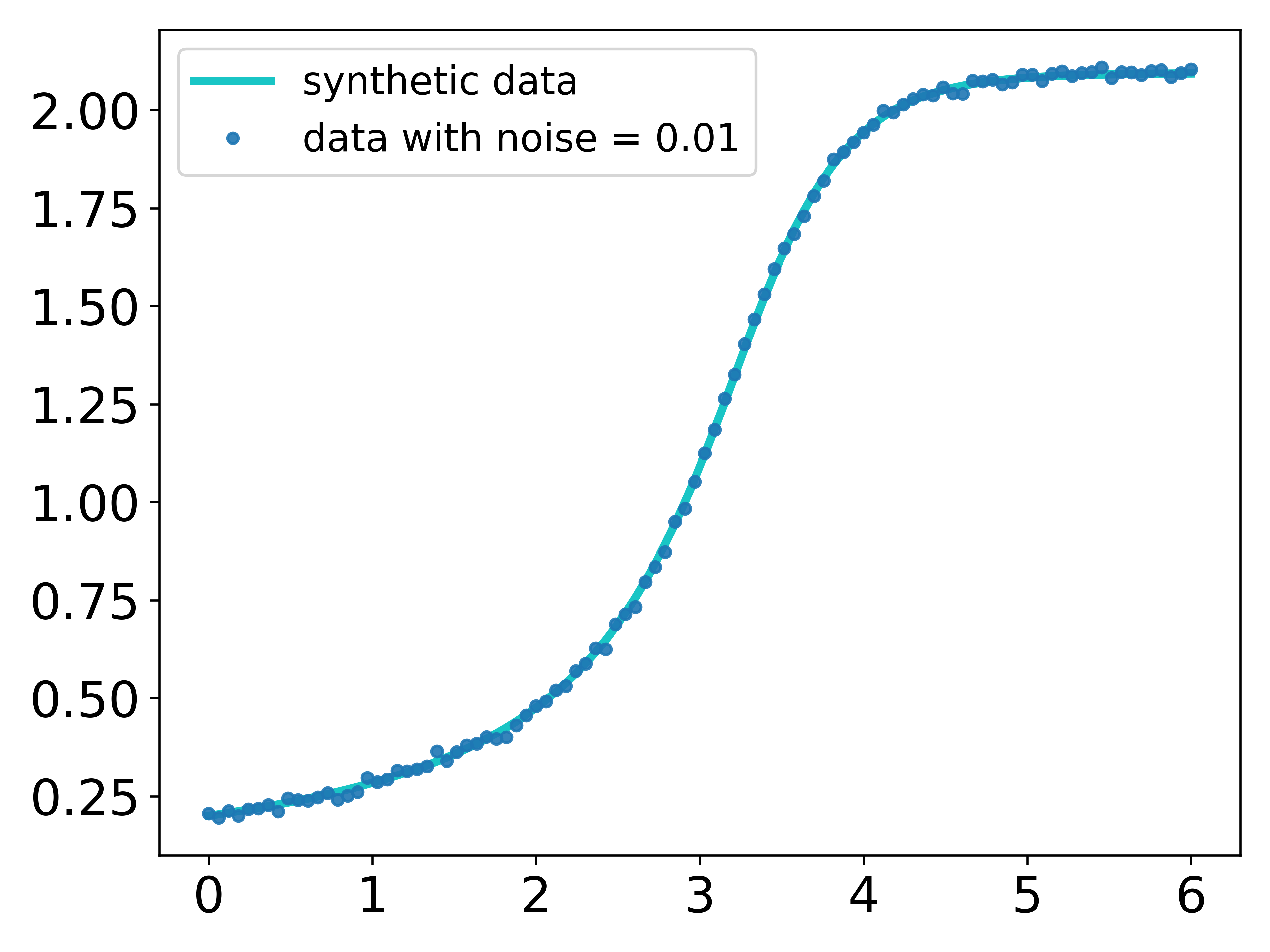} &
    \includegraphics[scale=0.29]{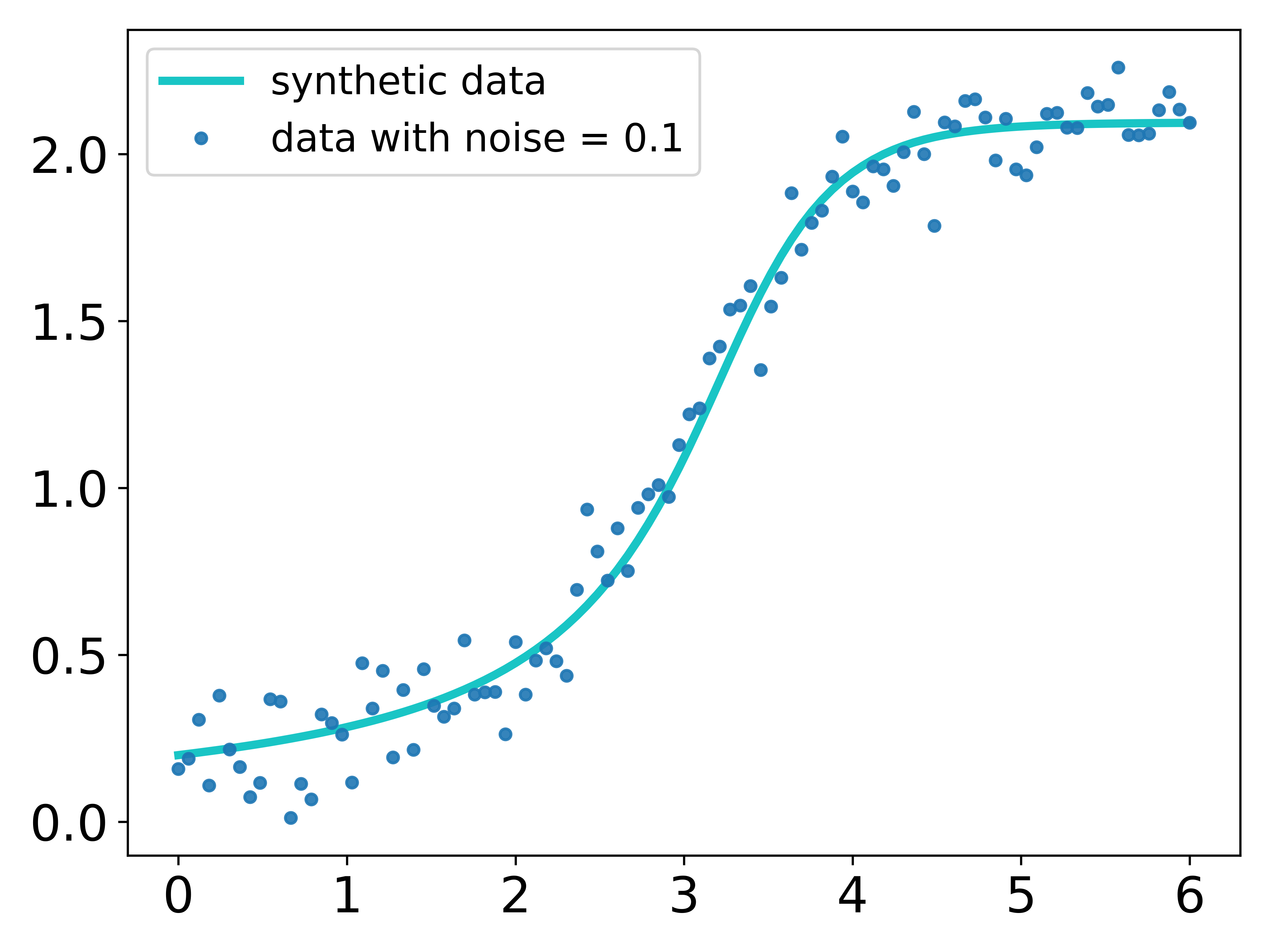} &
    \includegraphics[scale=0.29]{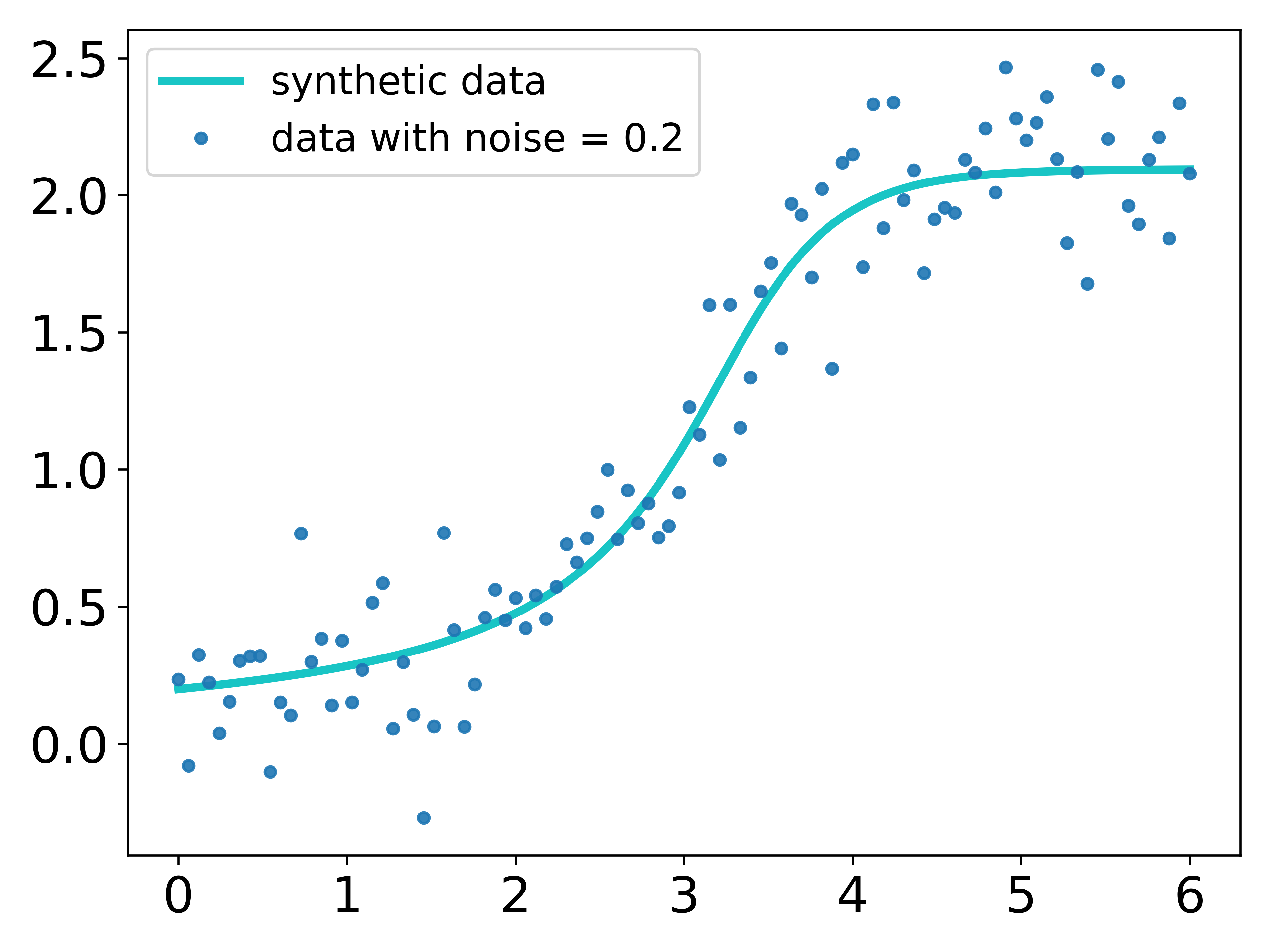} \\
    \includegraphics[scale=0.29]{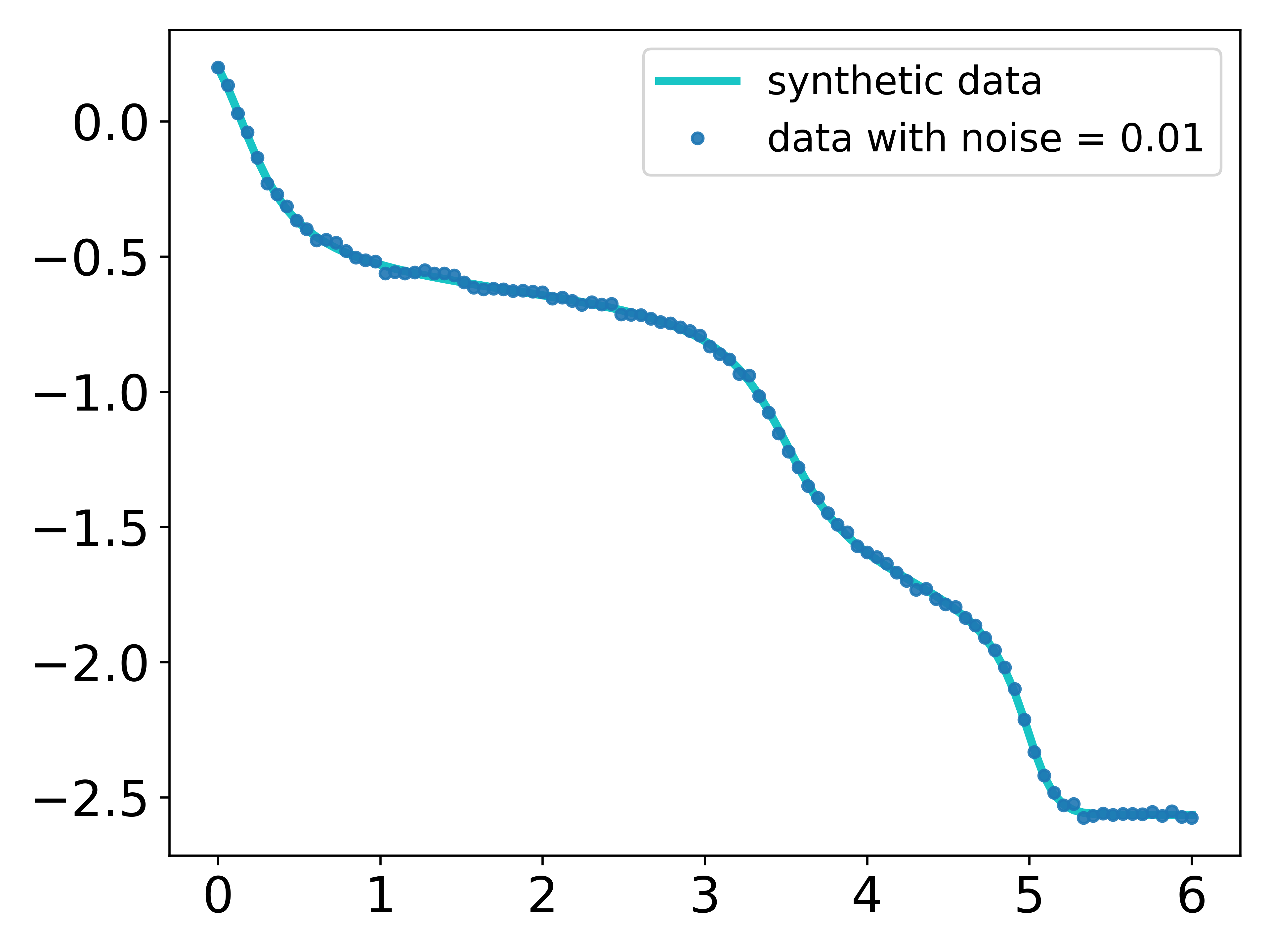} &
    \includegraphics[scale=0.29]{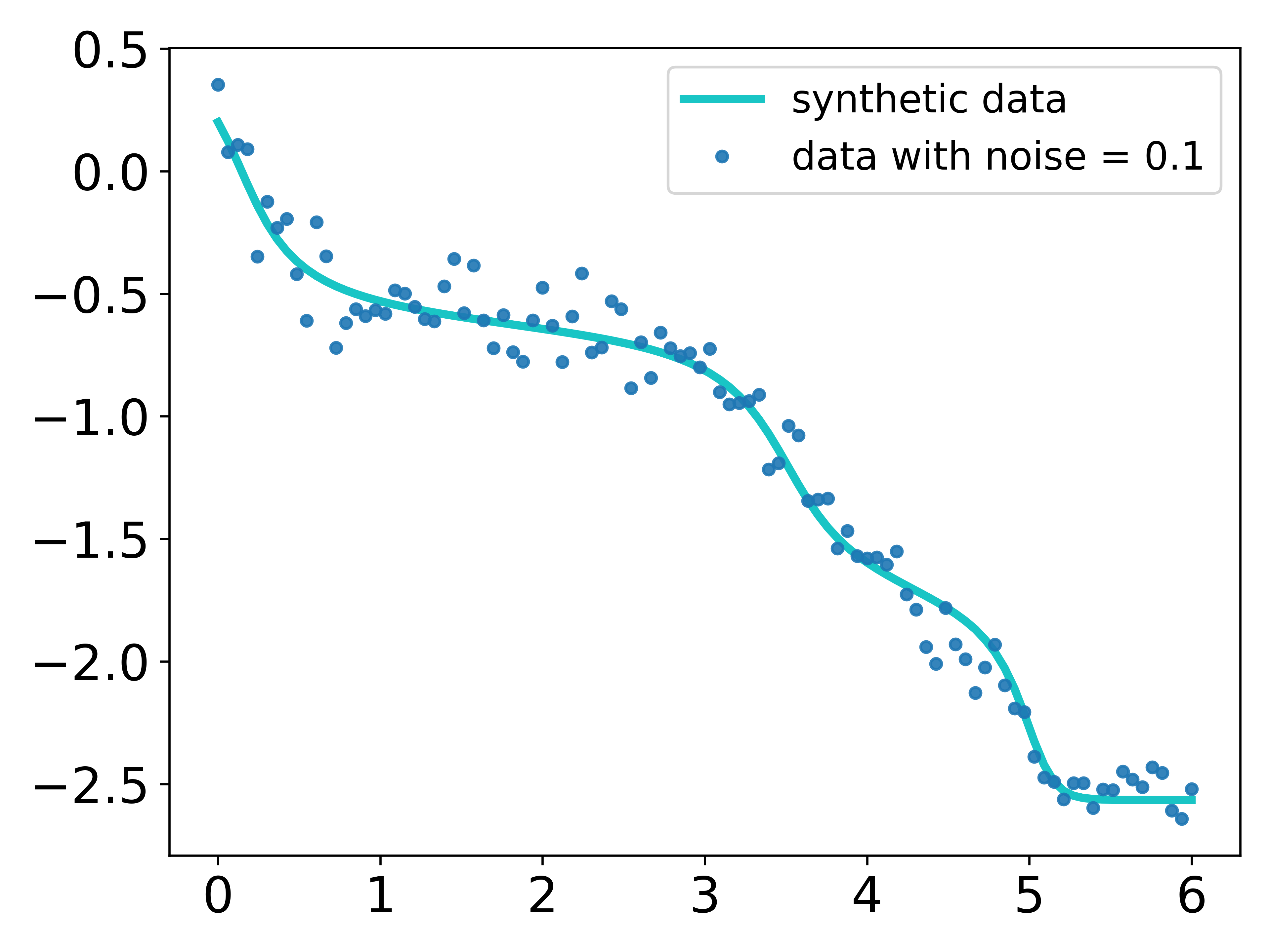} &
    \includegraphics[scale=0.29]{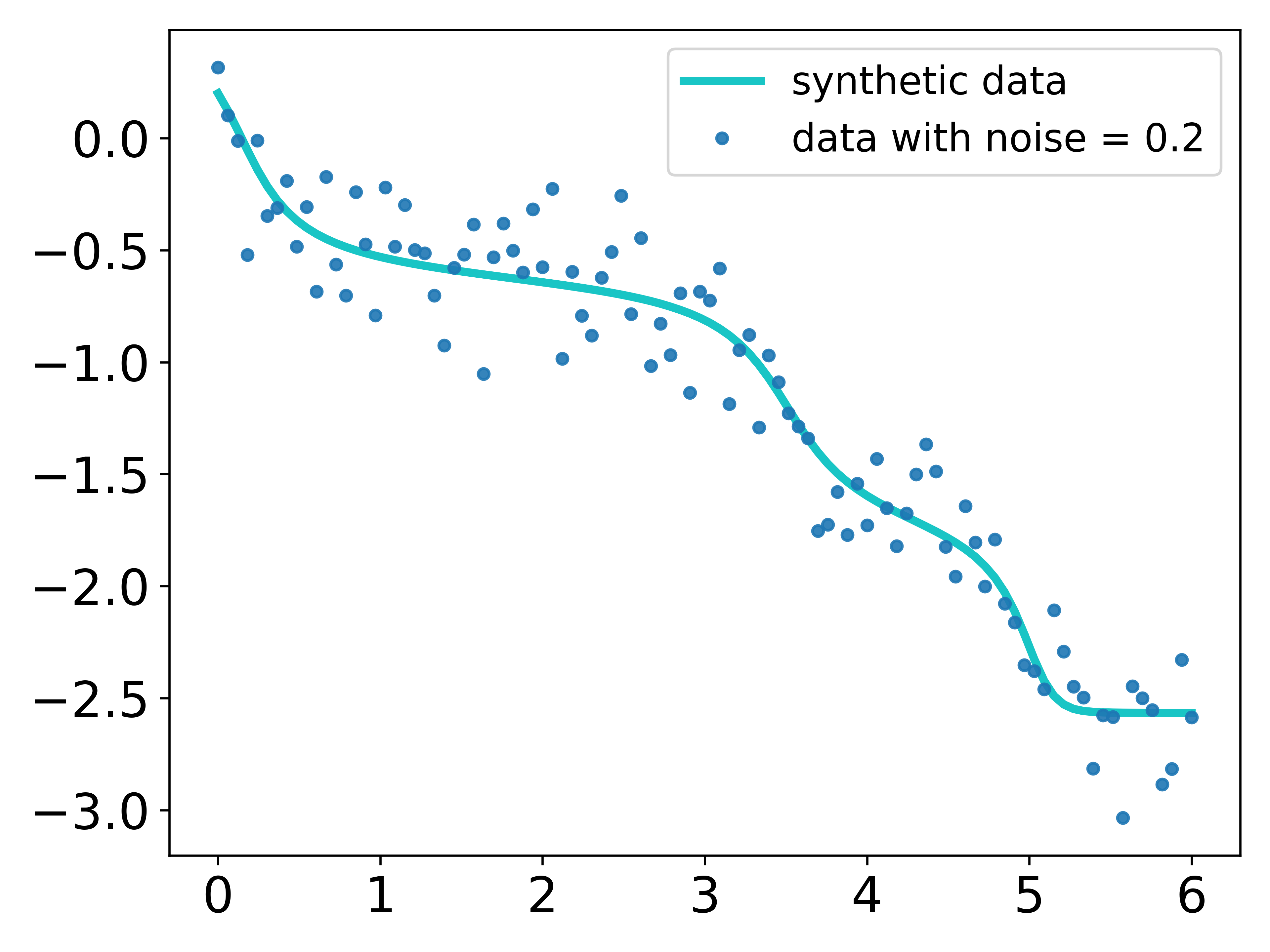} \\
    \end{tabular}
    \caption{$\textbf{m}_1$ (first row), $\textbf{m}_2$ (second row) synthetic data (light blue line) and noisy discrete synthetic data (blue dots) with different amount of noise: from left to right $\boldsymbol{\sigma} = [0.01, 0.1, 0.2]$. $x-$axis shows time. $y-$axis shows the values of the state variable.}
    \label{fig:synthetic-data-b1}
\end{figure}

Results obtained are summarized in Table \ref{tab:b1-errors} and illustrated via violin plots and regularization paths in Figures \ref{fig:b1-violin-plots} and \ref{fig:b1-regularization-paths} .

\begin{table}[H]\caption{Parameters and solutions mean relative error table. Each entry reports the mean relative error $\pm$ standard deviation on parameters and solutions. Each row is a different ground truth, each column is a different level of noise $\boldsymbol{\sigma} = [0.01, 0.1, 0.2]$}\label{tab:b1-errors}
\centering
\begin{tabular}{@{}lccc@{}}
\hline
\multicolumn{4}{c}{\textbf{Noise levels}} \\
\hline
M & Low & Medium & High \\
\hline
$\textbf{m}_1$ & 0.01 $\pm$ 0.00 & 0.02 $\pm$ 0.02 & 0.02 $\pm$ 0.01 \\
$\textbf{m}_2$ & 0.11 $\pm$ 0.00 & 0.12 $\pm$ 0.02 & 0.11 $\pm$ 0.05 \\
\hline
$\textbf{u}_1$ & 0.01 $\pm$ 0.00 & 0.05 $\pm$ 0.02 & 0.04 $\pm$ 0.09 \\
$\textbf{u}_2$ & 0.03 $\pm$ 0.00 & 0.05 $\pm$ 0.01 & 0.04 $\pm$ 0.16 \\
\hline
\end{tabular}
\end{table}

\begin{figure}[H]
\centering
\begin{tabular}{c c}
    \includegraphics[scale=0.25]{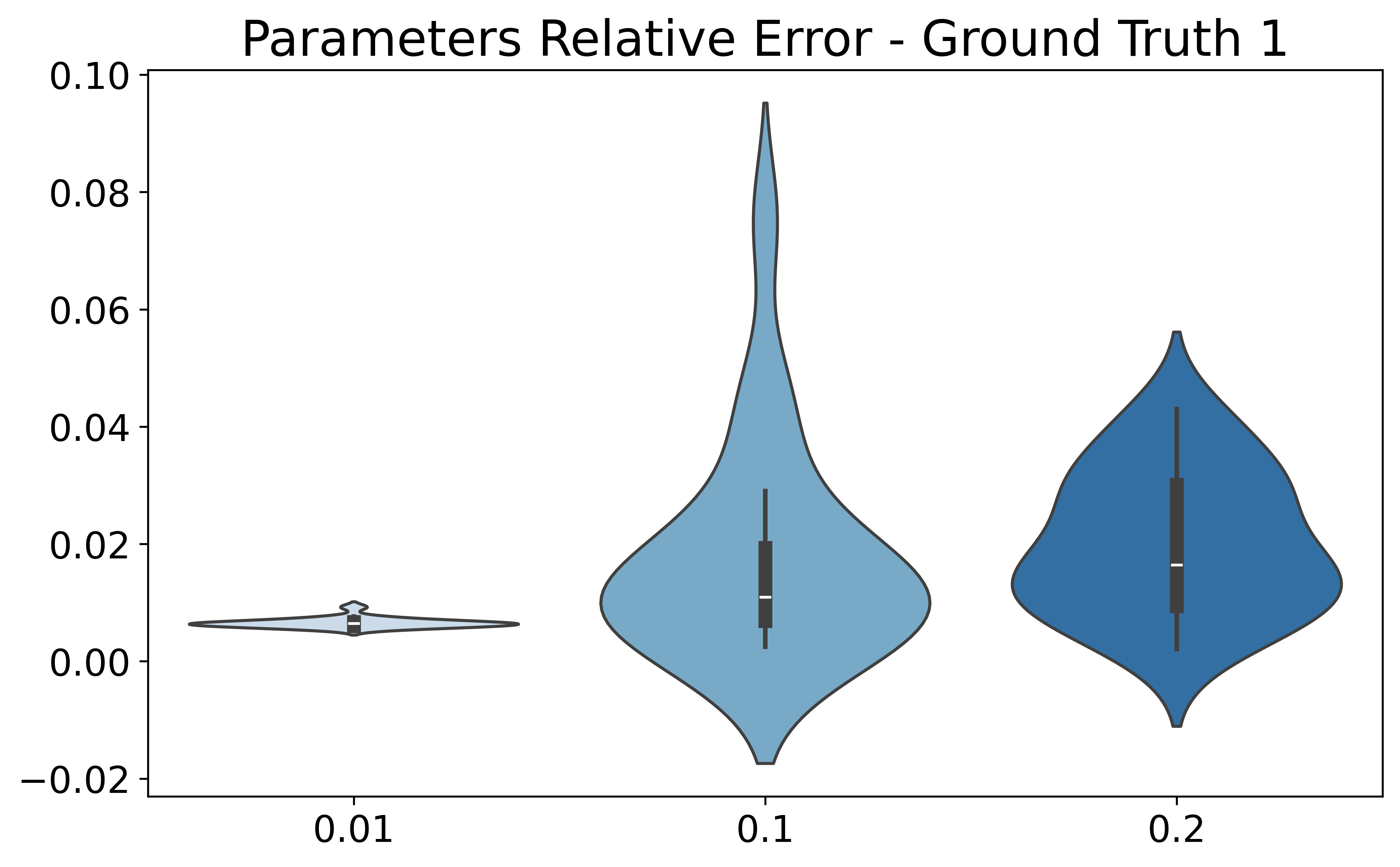} &
    \includegraphics[scale=0.25]{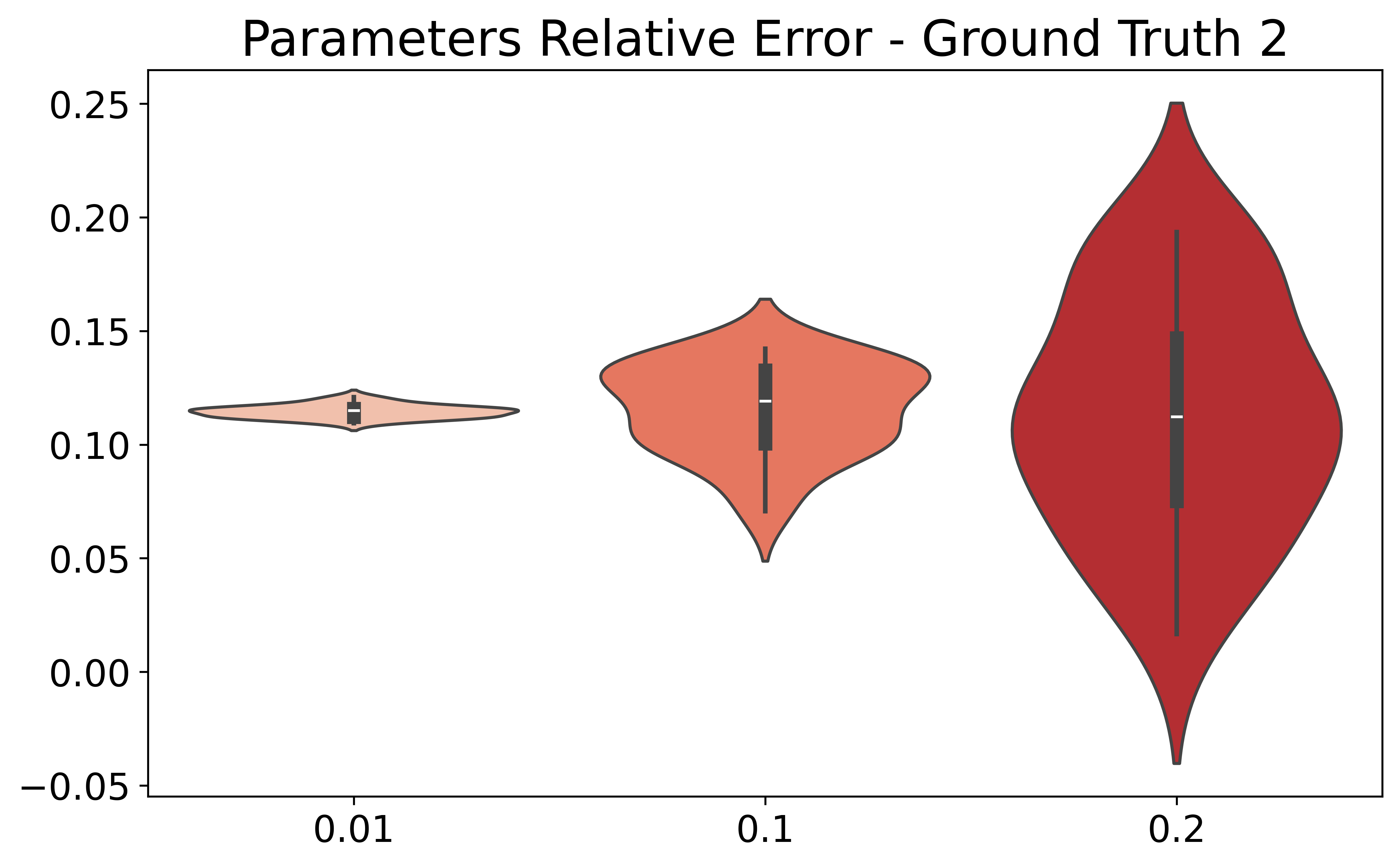} \\
    \includegraphics[scale=0.25]{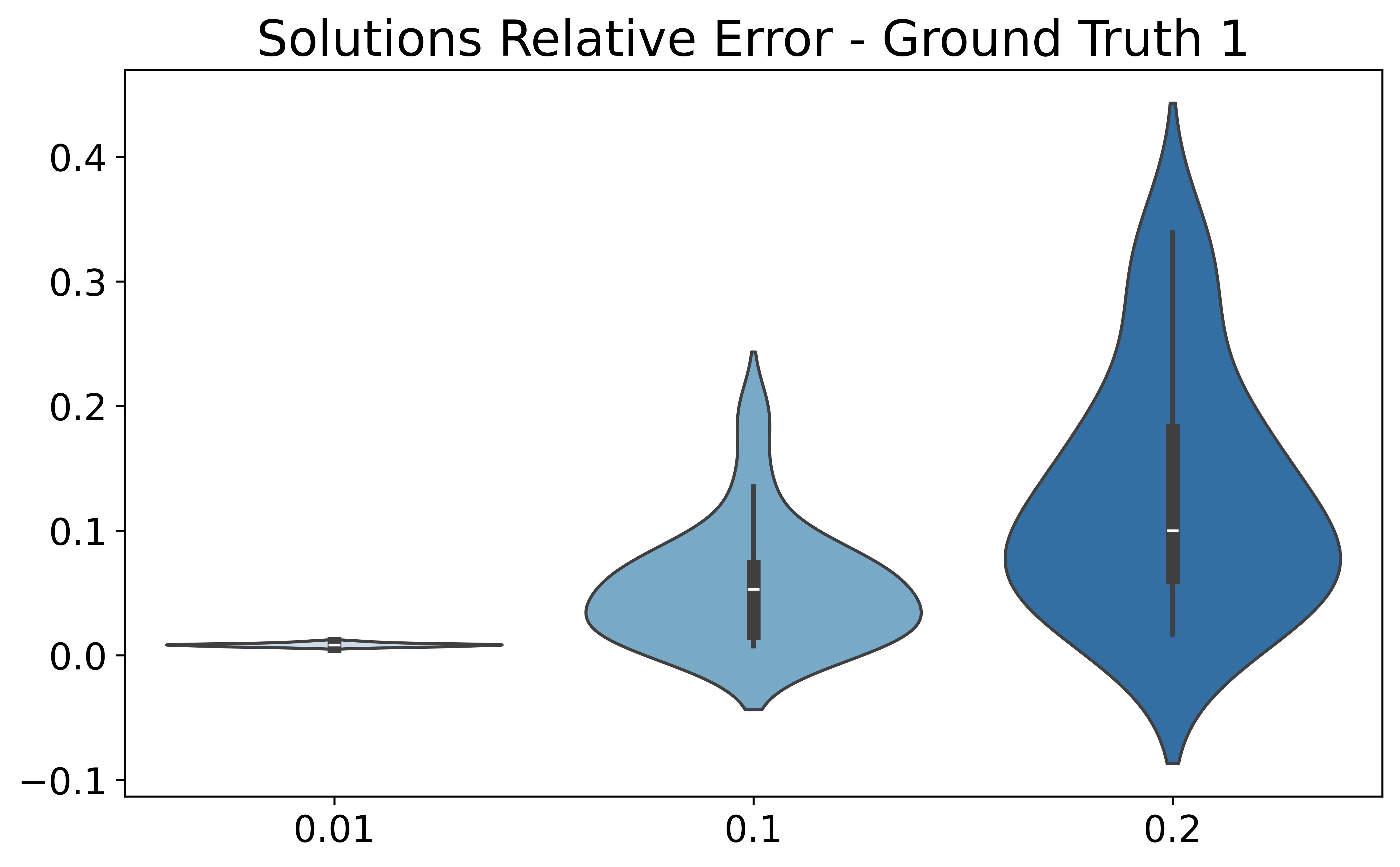} &
    \includegraphics[scale=0.25]{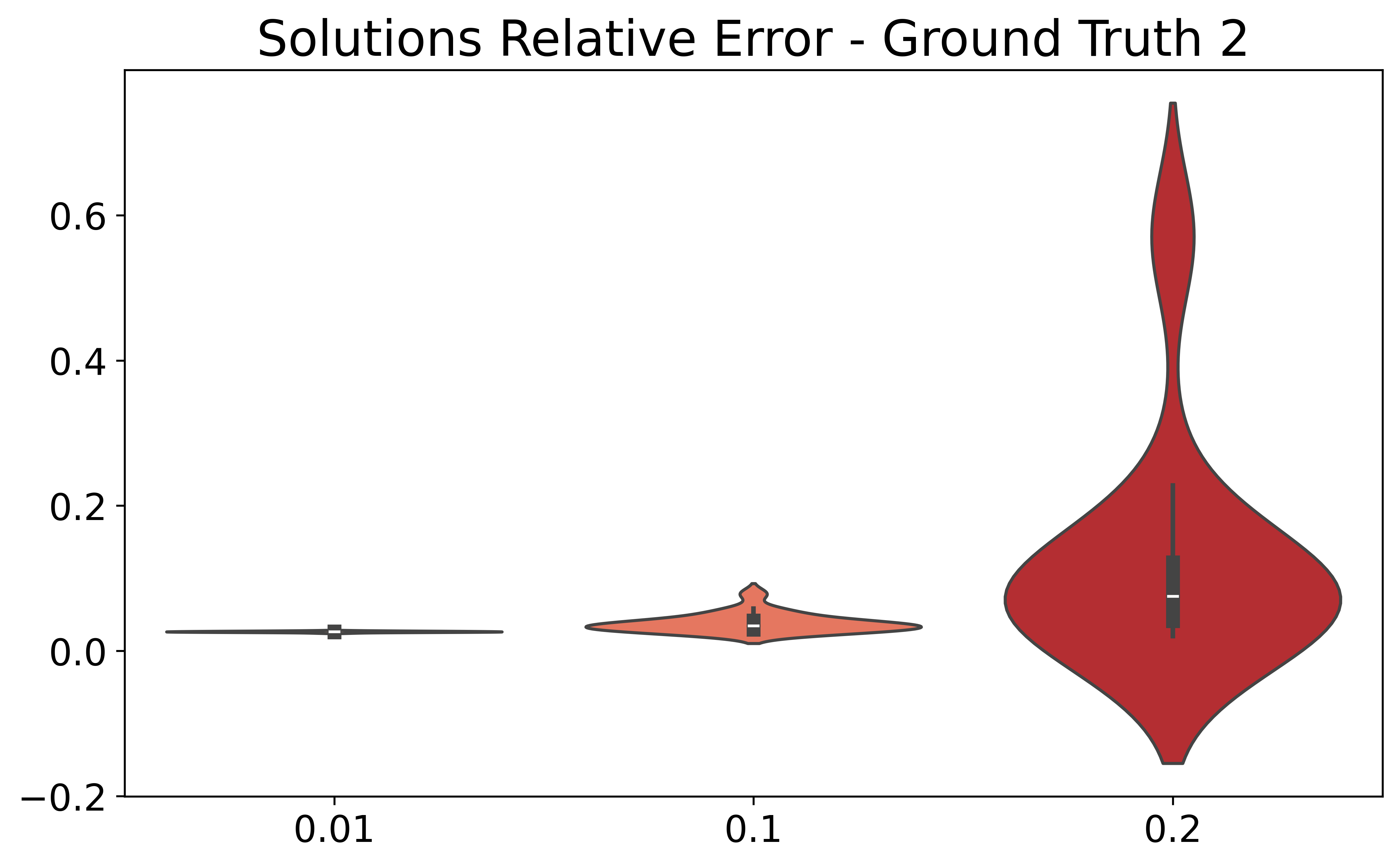}
    \end{tabular}
    \caption{Violin plots of relative errors on parameters and solutions computed with best regularization parameter $\alpha_q^*$ for each trial $q=1,...,n$. First row shows results about parameters $\textbf{m}$. Second row shows results about solutions $\textbf{u}$. $x$-axis reports different levels of noise $\boldsymbol{\sigma} = [0.01, 0.1, 0.2]$. $y$-axis reports the relative error.}
    \label{fig:b1-violin-plots}
\end{figure}

\begin{figure}[H]
\centering
\begin{tabular}{c c c}
    \includegraphics[scale=0.3]{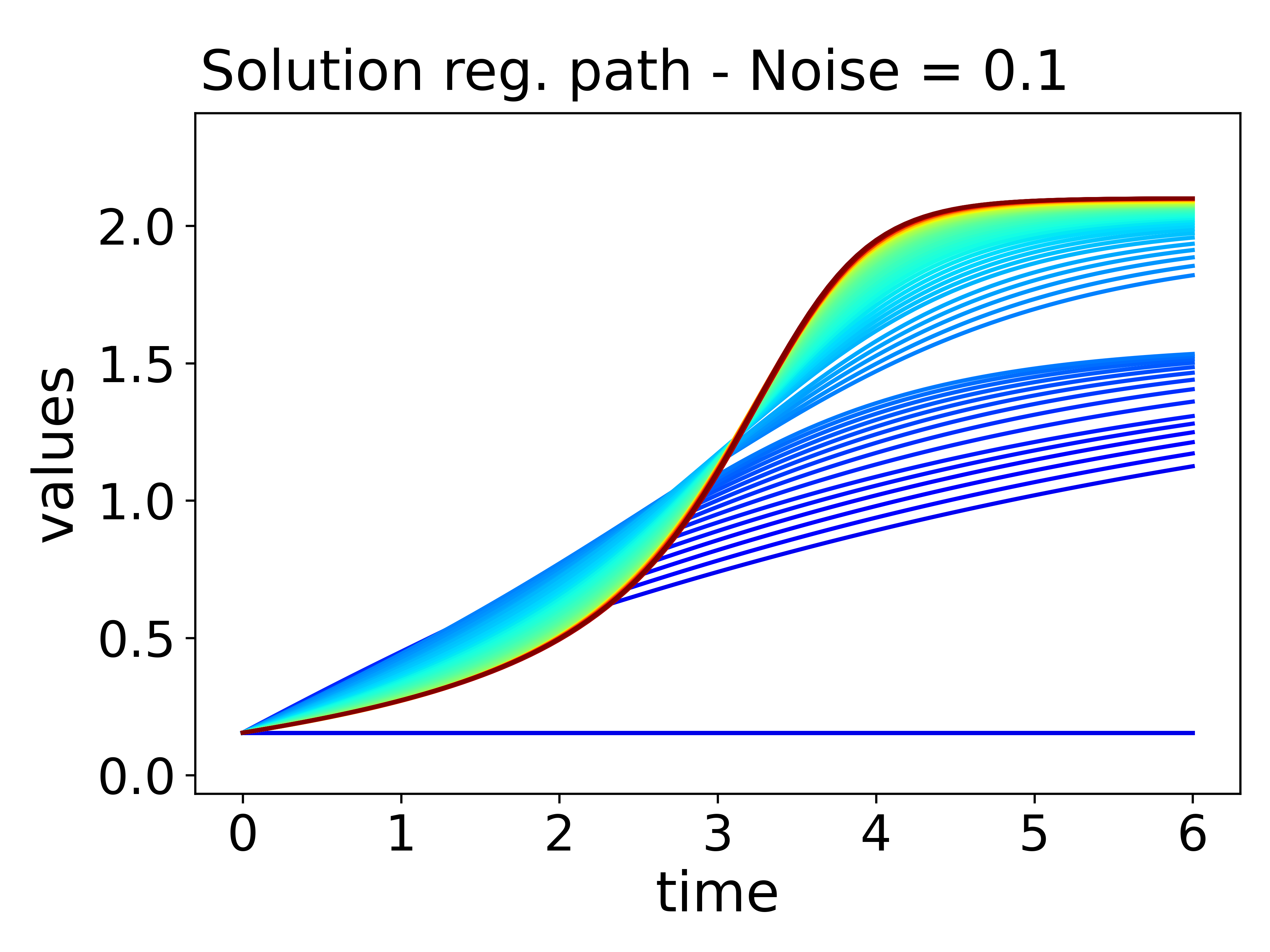} &
    \includegraphics[scale=0.3]{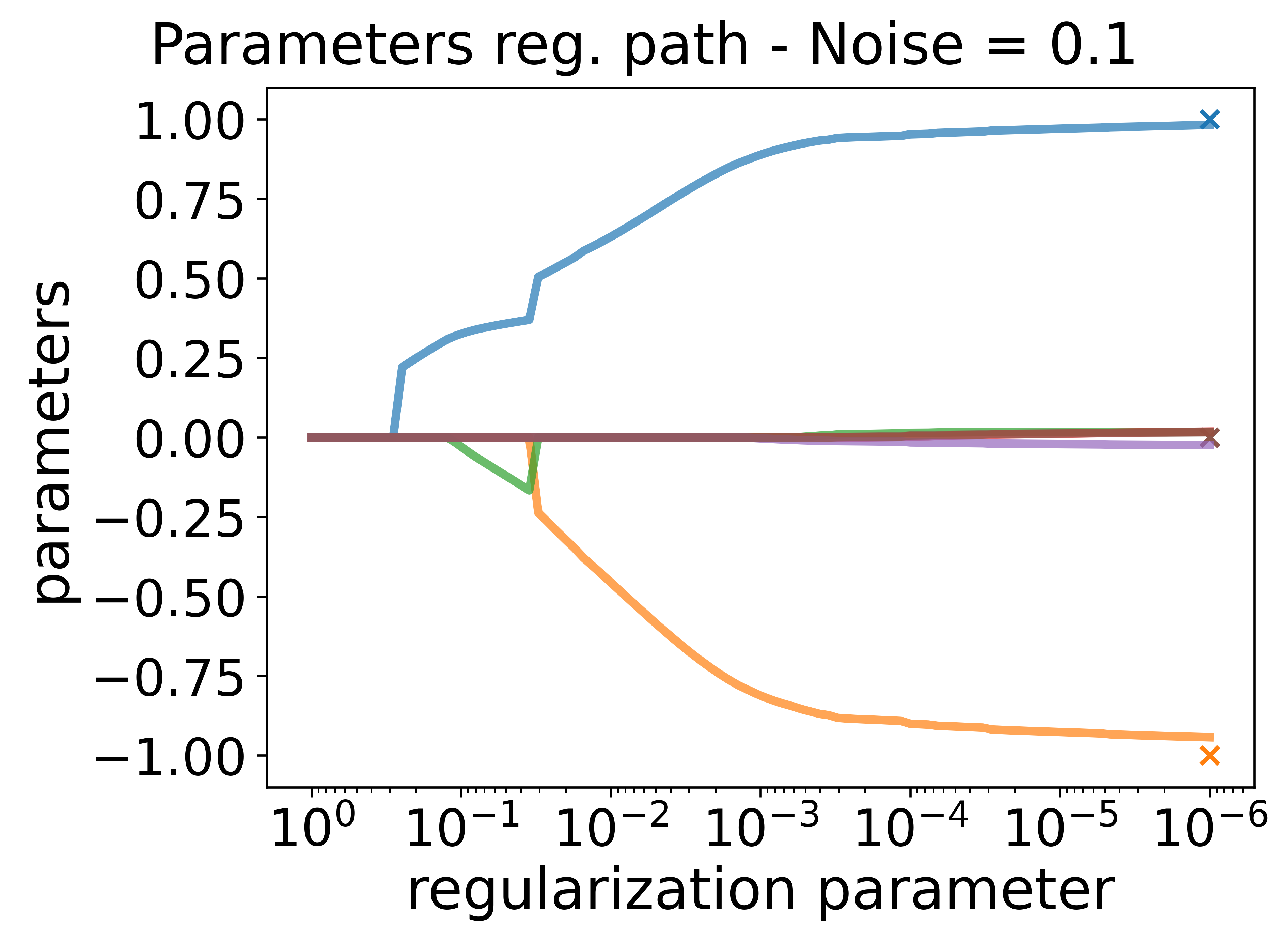} & 
    \includegraphics[scale=0.3]{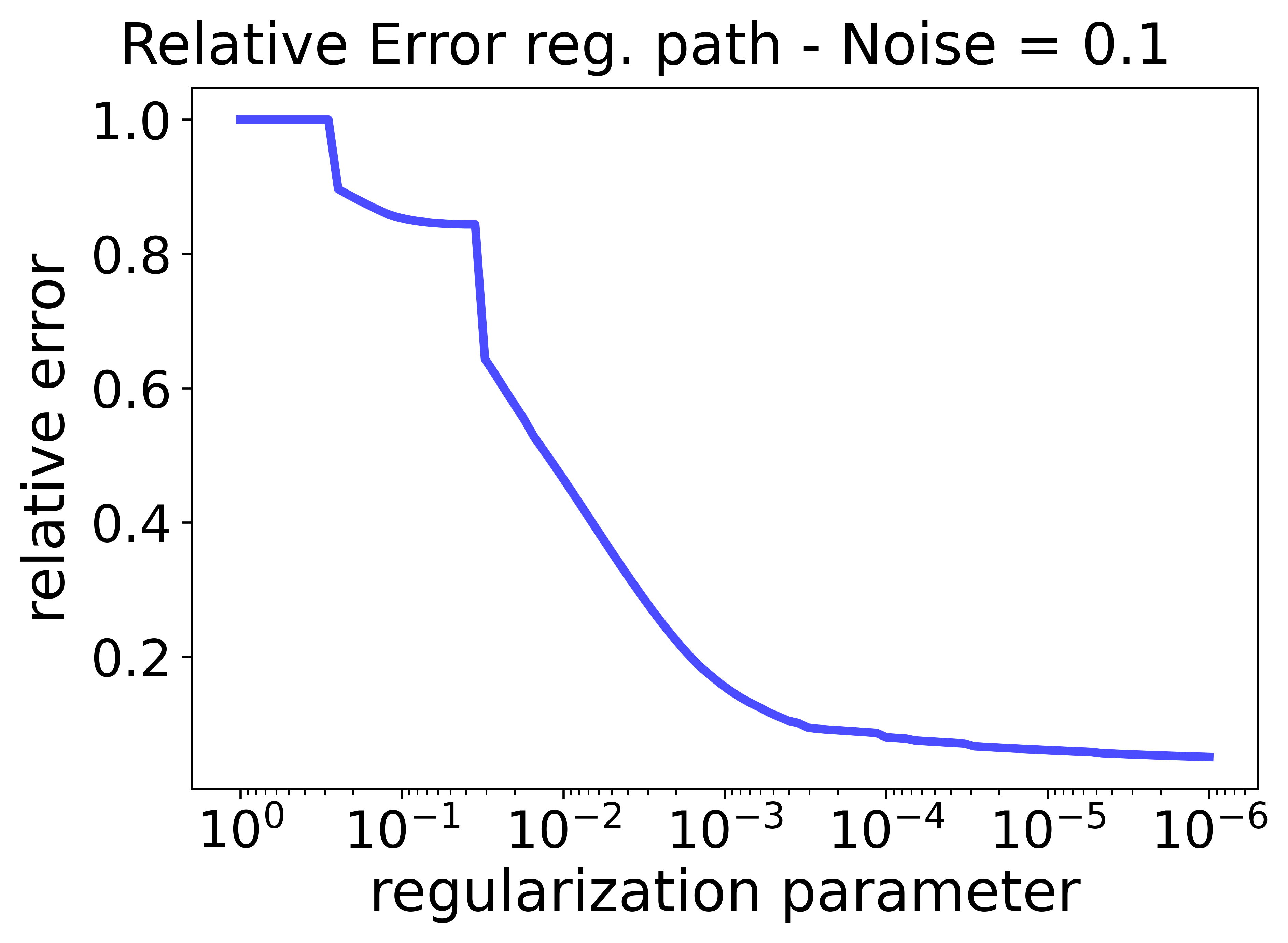} \\
    \includegraphics[scale=0.3]{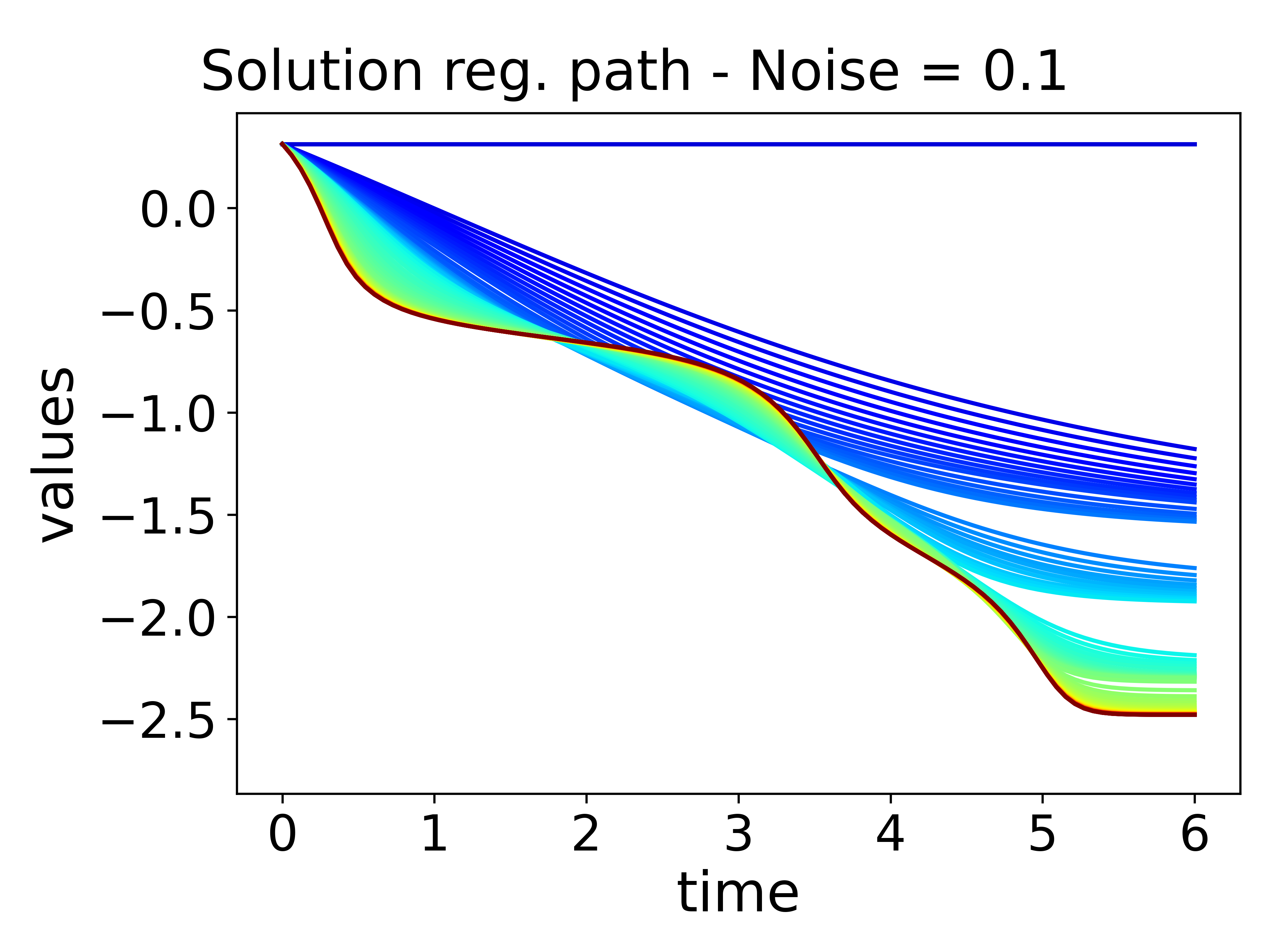} &
    \includegraphics[scale=0.3]{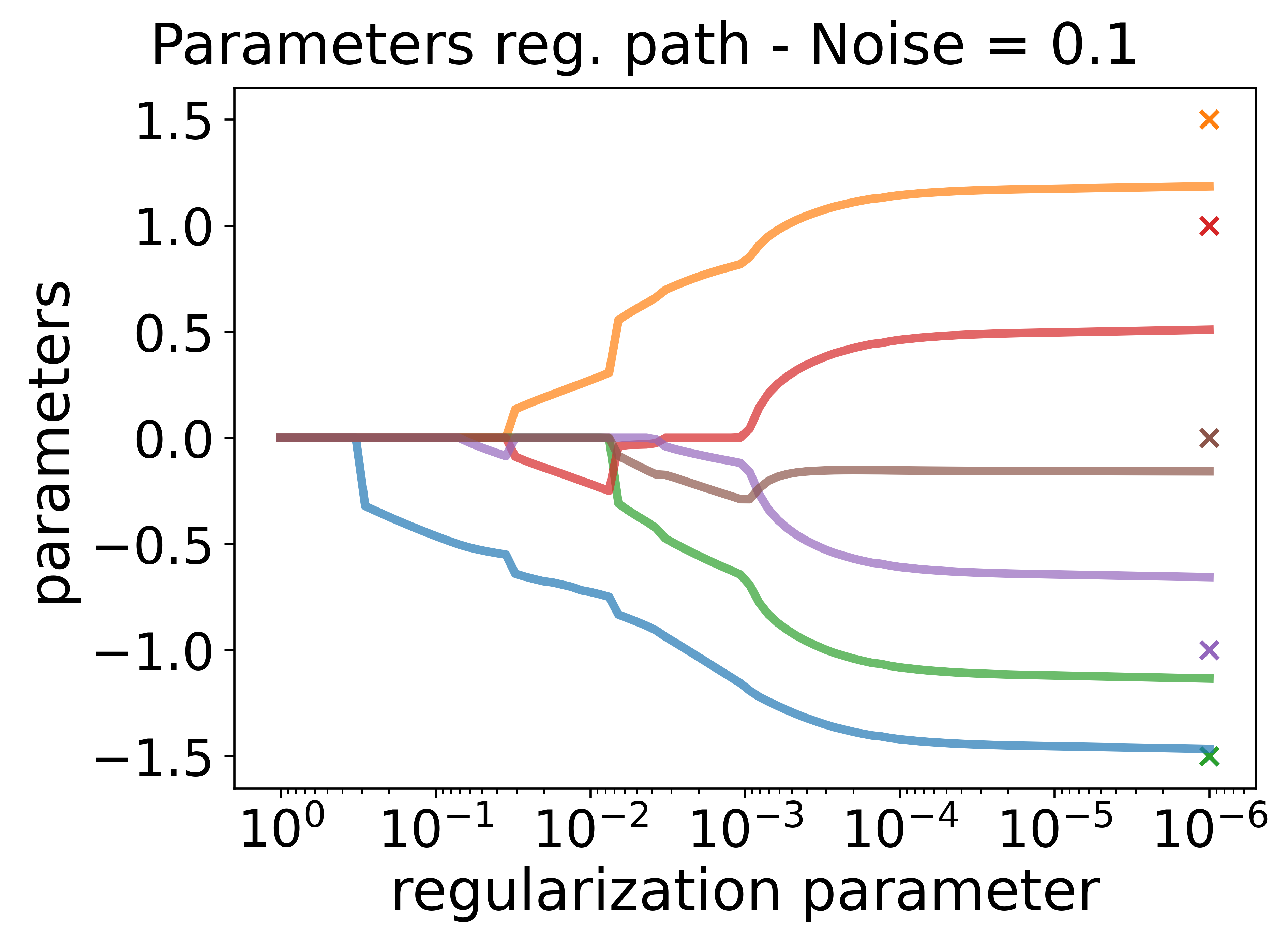} & 
    \includegraphics[scale=0.3]{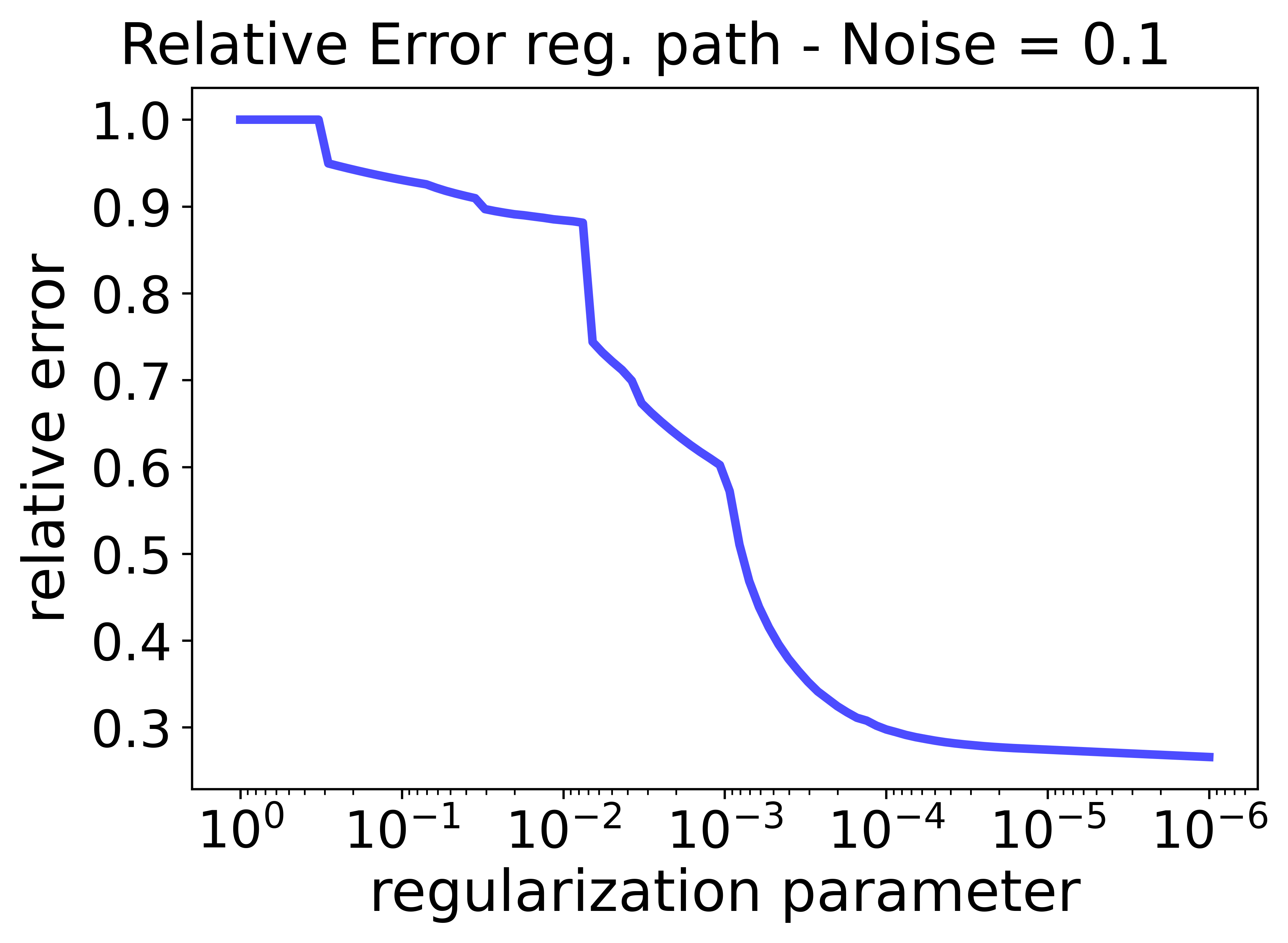}
    \end{tabular}
    \caption{Regularization paths with level of noise $\sigma = 0.1$. First row shows results for $\textbf{m}_1$ \eqref{b1-m1}, second row shows results for $\textbf{m}_2$ \eqref{b1-m2}. The first column reports the regularization paths of the solutions. The colors of the curves transition from blue ($l=0$) to red ($l=99$). The $x-$axis reports the time, the $y-$axis the value of the state variable. The second column reports the regularization paths of the parameters. $x-$axis shows the regularization parameter, $y-$axis shows the parameter values. The ground truth values are symbolized by a cross at the last regularization parameter; the curves and crosses of the same color correspond to the same parameter. The third column shows the regularization path of the relative error on parameters. $x-$axis shows the regularization parameter, $y-$axis shows the relative error.}
    \label{fig:b1-regularization-paths}
\end{figure}

From Table \ref{tab:b1-errors} and the violin plots shown in Figure \ref{fig:b1-violin-plots}, it can be observed that the proposed method demonstrates improved performance when the ground truth is sparse ($\textbf{m}_1$), with relative errors approximately one order of magnitude lower than in the less sparse case ($\textbf{m}_2$). The approach also shows robustness to noise, as reconstruction accuracy remains stable across different noise levels. Despite variations in parameter recovery, the relative error on the reconstructed solutions remains low overall. Regularization paths are smooth and closely match the synthetic data, indicating numerical stability. These results support the self-consistency of the method, which is specifically designed to promote sparse solutions through the structure of the regularization term.

Figure \ref{fig:semi-convergence} is obtained by considering a noise level of $\sigma = 0.2$, a choice that allows us to operate in a context where regularization is indispensable, while keeping the perturbed data sufficiently close to the exact data. The figure shows a representative example of the semi-convergence phenomenon, which often occurs in linear inverse problems when there is numerical instability due to noise.

\begin{figure}[H]
\centering
\begin{tabular}{c c}
    \includegraphics[scale=0.3]{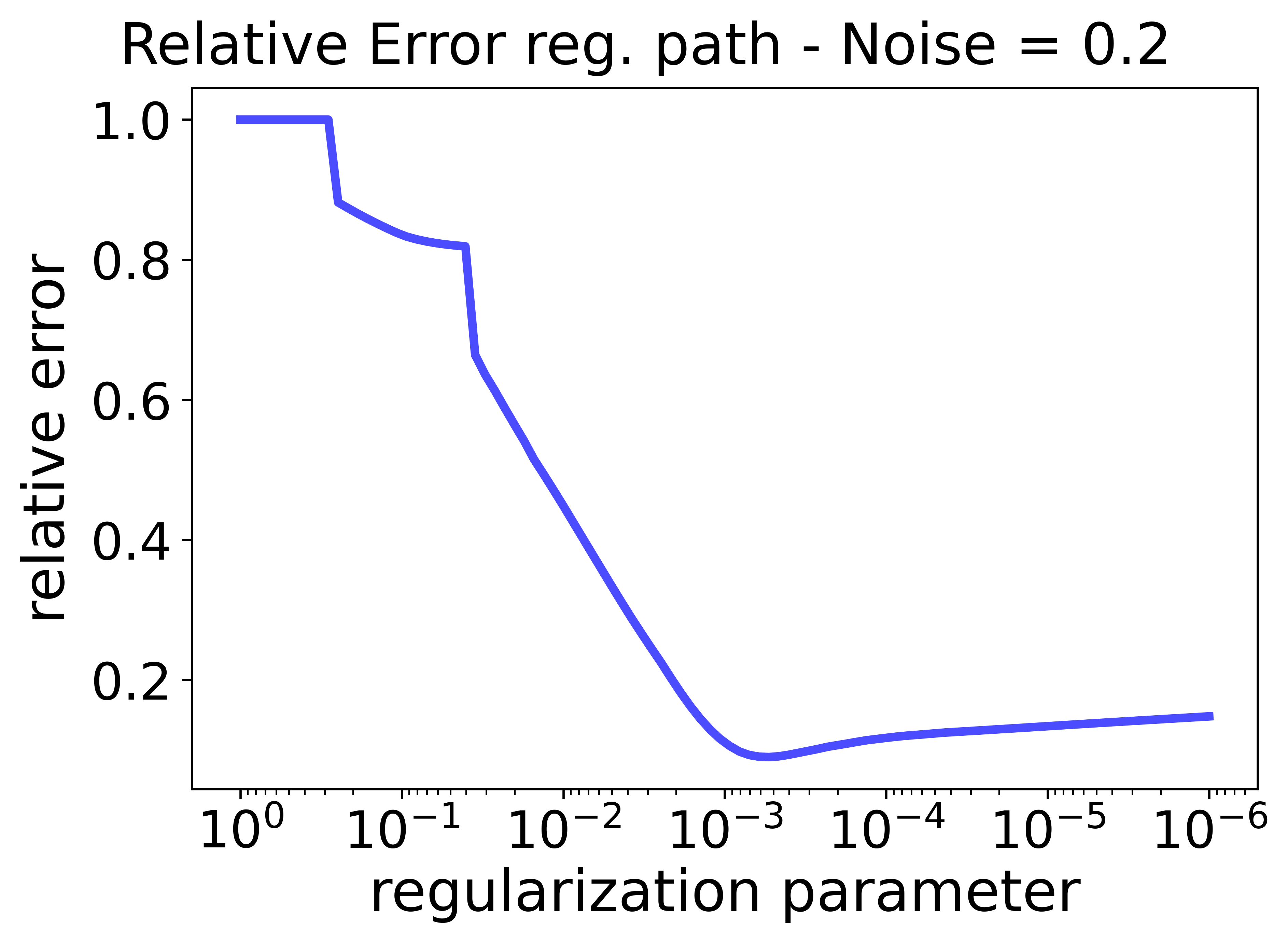} &
    \includegraphics[scale=0.3]{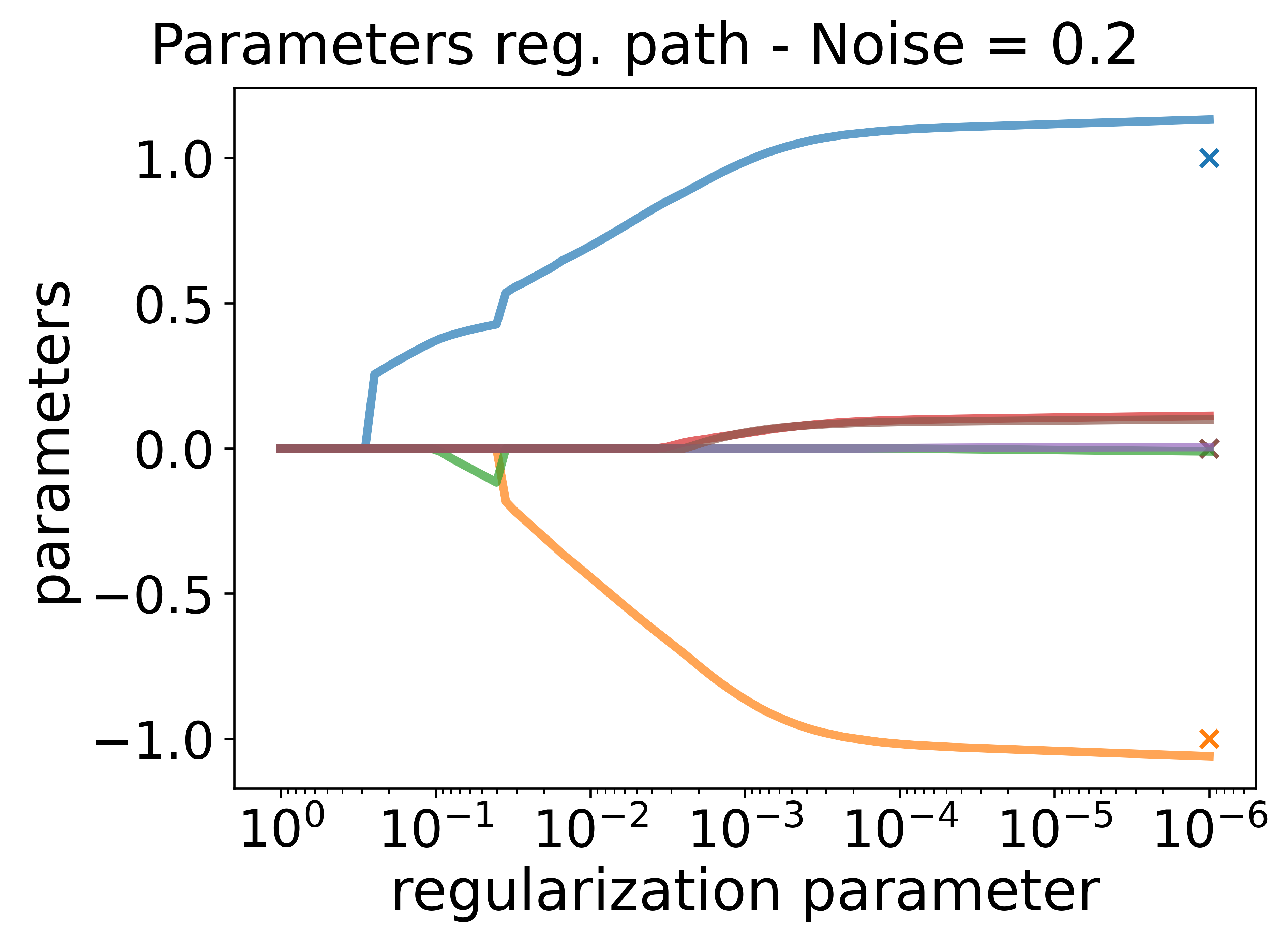} \\
    \includegraphics[scale=0.3]{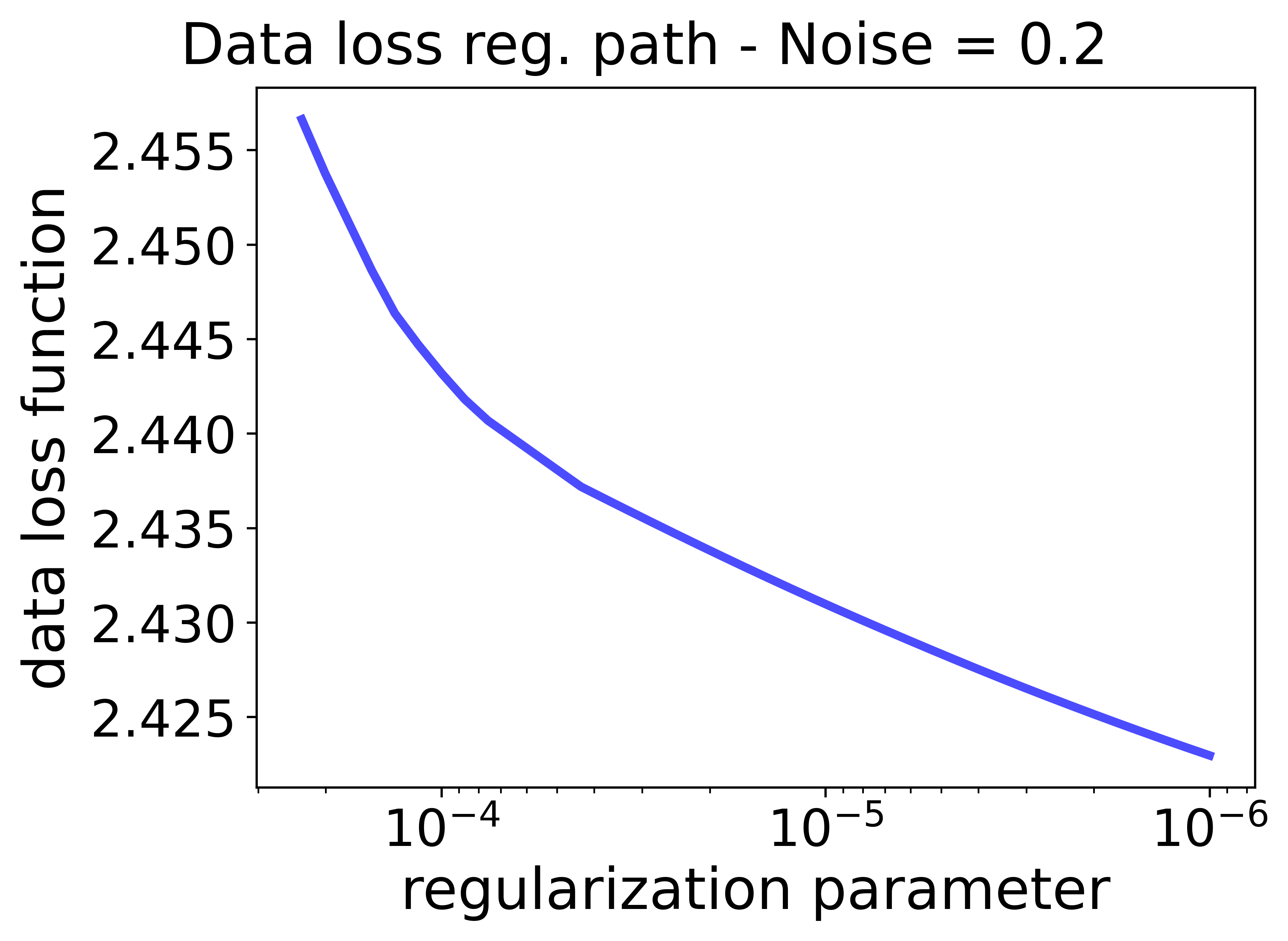} &
    \includegraphics[scale=0.3]{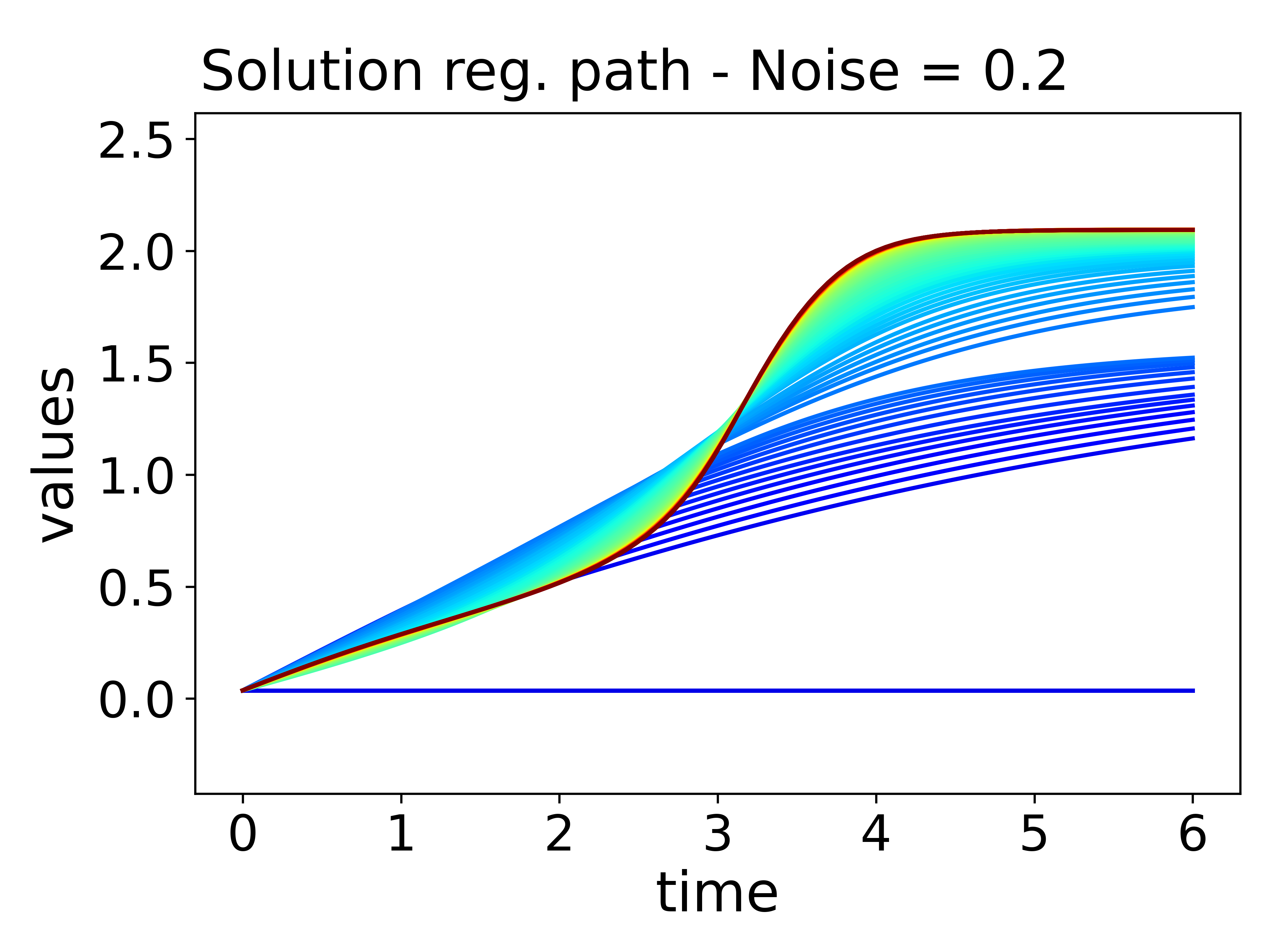}
    \end{tabular}
    \caption{Regularization paths with level of noise $\sigma = 0.2$ for $\textbf{m}_1$ \eqref{b1-m1}. Top left: regularization path of the relative error on parameters. $x-$axis shows the regularization parameter, $y-$axis shows the relative error. Top right: regularization paths of the parameters. $x-$axis shows the regularization parameter, $y-$axis shows the parameter values. The ground truth values are symbolized by a cross at the last regularization parameter; the curves and crosses of the same color correspond to the same parameter. 
    Bottom left: terminal part of regularization path of the data loss function \eqref{data loss}.
    Bottom right: regularization paths of the solutions. The colors of the curves transition from blue ($l=0$) to red ($l=99$). The $x-$axis reports the time, the $y-$axis the value of the state variable.}
    \label{fig:semi-convergence}
\end{figure}

In particular, the top left panel \ref{fig:semi-convergence}, shows how the relative error on the parameters initially decreases as the regularization parameter decreases. It then reaches a minimum and increases again when the parameter becomes too small. At this stage, the algorithm begins to reconstruct the spurious noise components.
The upper-right panel represents the evolution of the parameters along the regularization path and confirms this semi-convergent behavior. The optimal configuration of the parameters, which is the most consistent with the ground truth, is not obtained at extreme values of the regularization parameter. Rather, it is located in an intermediate region of the path.
Despite this, the algorithm continues to optimize the loss function, as shown in the lower left (scaled and zoomed in to the last part of the regularization path) panel. The bottom right panel shows the evolution of the forward solutions of the Cauchy problem \eqref{Cauchy problem} along the regularization path. The gradual convergence to the synthetic data curve, shown in Figure \ref{fig:synthetic-data-b1}, again highlights the ill-posed nature of the problem.

\section{Conclusions}
In this paper, we proposed a homotopy-based regularization path method for solving implicit inverse problems. Specifically, we aimed to determine a nonlinear functional of physical quantities directly from their noisy measurements. Assuming a parametric form for the functional \eqref{State equation}, the implicit inverse problem was reformulated as a parameter optimization problem.
By introducing a regularized objective function \eqref{Objective Functional}, the problem was further recast as a discretized, regularized, and constrained minimization problem \eqref{Minimum regluarized constrained problem discrete}.
Fixing the regularization parameter, we employed a variational approach and derived a gradient descent-based algorithm \eqref{Gradient Step}. To efficiently compute large Jacobians \eqref{lagrangian p d m}, we used the adjoint state method \eqref{adjoint State System}. To approach the solution of the non-regularized constrained problem \eqref{Minimum constrained problem discrete}, we applied a continuation strategy based on the regularization path and warm-start initialization \eqref{regularization parameter}, culminating in the proposed algorithm summarized in pseudocode (\ref{alg:RASM}).

The method was then specialized for latent dynamics discovery, where the unknown relationships were modeled as systems of ordinary differential equations (ODEs), with time series as the observed data. The unknown nonlinear functional was represented as a linear combination of nonlinear basis functions \eqref{basis function}, as shown in \eqref{Derivative function sparse}, and the forward problem was cast as a Cauchy problem \eqref{Cauchy problem}.
Numerical experiments in a synthetic one-dimensional setting validated the effectiveness of the proposed method. In tests involving periodic functions, the method showed strong performance, especially when the ground truth was scattered. These results confirm the internal consistency of the approach, in line with the sparsity-promoting regularization employed. The experiments also underscored the importance of regularization for problem stability, particularly in the presence of semi-convergence phenomena.
%The results also show the ability of the method to find the presence of bifurcations in the regularization path, allowing multiple branches to be tracked.

A key limitation of the method is the design of the candidate basis function library, which must be expressive enough to capture the system dynamics without introducing unnecessary complexity. Additionally, the current implementation requires manual computation of the basis function derivatives, which becomes burdensome as the basis set grows.

Future work will focus on improving the algorithm's efficiency, particularly in the context of expanding the basis function library. This includes leveraging automatic differentiation to simplify derivative computations and developing systematic methods for scaling the library. We also plan to extend the methodology to more complex dynamic systems, including those involving multiple physical quantities, spatiotemporal phenomena governed by partial differential equations, and systems subject to inequality constraints.

\section*{References}

\bibliographystyle{iopart-num}%apa{plain}%{iopart-num}
% \bibliographystyle{unsrt}

%\bibliography{mybiblio}

\begin{thebibliography}{10}
\expandafter\ifx\csname url\endcsname\relax
  \def\url#1{{\tt #1}}\fi
\expandafter\ifx\csname urlprefix\endcsname\relax\def\urlprefix{URL }\fi
\providecommand{\eprint}[2][]{\url{#2}}

\bibitem{mosegaard1999probabilistic}
Mosegaard K and Rygaard-Hjalsted C 1999 {\em Inverse problems\/} {\bf 15} 573

\bibitem{nino2018local}
Nino-Ruiz E~D, Ardila C and Capacho R 2018 {\em Soft Computing\/} {\bf 22} 4819--4832

\bibitem{brunton2016discovering}
Brunton S~L, Proctor J~L and Kutz J~N 2016 {\em Proceedings of the national academy of sciences\/} {\bf 113} 3932--3937

\bibitem{chen2021physics}
Chen Z, Liu Y and Sun H 2021 {\em Nature communications\/} {\bf 12} 6136

\bibitem{heinonen2018learning}
Heinonen M, Yildiz C, Mannerstr{\"o}m H, Intosalmi J and L{\"a}hdesm{\"a}ki H 2018 Learning unknown ode models with gaussian processes {\em International conference on machine learning\/} (PMLR) pp 1959--1968

\bibitem{lorenzi2018constraining}
Lorenzi M and Filippone M 2018 Constraining the dynamics of deep probabilistic models {\em International Conference on Machine Learning\/} (PMLR) pp 3227--3236

\bibitem{dondelinger2013ode}
Dondelinger F, Husmeier D, Rogers S and Filippone M 2013 Ode parameter inference using adaptive gradient matching with gaussian processes {\em Artificial intelligence and statistics\/} (PMLR) pp 216--228

\bibitem{williams2006gaussian}
Williams C~K and Rasmussen C~E 2006 {\em Gaussian processes for machine learning\/} vol~2 (MIT press Cambridge, MA)

\bibitem{rielly2025mock}
Rielly V, Lahouel K, Lew E, Fisher N, Haney V, Wells M and Jedynak B in press {\em stat\/} {\bf 1050} 8

\bibitem{plessix2006review}
Plessix R~E 2006 {\em Geophysical Journal International\/} {\bf 167} 495--503

\bibitem{bergounioux2019position}
Bergounioux M, Bretin {\'E} and Privat Y 2019 {\em Inverse Problems\/} {\bf 35} 074003

\bibitem{alexanderian2021optimal}
Alexanderian A 2021 {\em Inverse Problems\/} {\bf 37} 043001

\bibitem{baayen2019overview}
Baayen J, Becker B, van Heeringen K~J, Miltenburg I, Piovesan T, Rauw J, den Toom M and VanderWees J 2019 {\em IFAC-PapersOnLine\/} {\bf 52} 73--80

\bibitem{keller2024ai}
Keller J~D and Potthast R 2024 {\em arXiv preprint arXiv:2406.00390\/}

\bibitem{burger2020data}
Burger M, Pietschmann J~F and Wolfram M~T 2020 {\em Inverse Problems\/} {\bf 36} 064003

\bibitem{rothermel2021solving}
Rothermel D and Schuster T 2021 {\em Inverse Problems\/} {\bf 37} 045014

\bibitem{migorski2019inverse}
Mig{\'o}rski S, Khan A~A and Zeng S 2019 {\em Inverse Problems\/} {\bf 35} 035004

\bibitem{watson1989modern}
Watson L~T and Haftka R~T 1989 {\em Computer Methods in Applied Mechanics and Engineering\/} {\bf 74} 289--305

\bibitem{Dedieu2015}
Dedieu J~P 2015 {\em Newton-Raphson Method\/} (Springer Berlin Heidelberg) pp 1023--1028 ISBN 978-3-540-70529-1 \urlprefix\url{https://doi.org/10.1007/978-3-540-70529-1-374}

\bibitem{friedman2010regularization}
Friedman J, Hastie T and Tibshirani R 2010 {\em Journal of statistical software\/} {\bf 33} 1

\bibitem{hadamard1902problemes}
Hadamard J 1902 {\em Princeton university bulletin\/}  49--52

\bibitem{tikhonov1977solutions}
Tikhonov A~N and Arsenin V 1977 {\em (No Title)\/}

\bibitem{rudin1992nonlinear}
Rudin L~I, Osher S and Fatemi E 1992 {\em Physica D: nonlinear phenomena\/} {\bf 60} 259--268

\bibitem{bertsekas2014constrained}
Bertsekas D~P 2014 {\em Constrained optimization and Lagrange multiplier methods\/} (Academic press)

\bibitem{rosasco2004loss}
Rosasco L, De~Vito E, Caponnetto A, Piana M and Verri A 2004 {\em Neural computation\/} {\bf 16} 1063--1076

\bibitem{hanke1995convergence}
Hanke M, Neubauer A and Scherzer O 1995 {\em Numerische Mathematik\/} {\bf 72} 21--37

\bibitem{scherzer1995convergence}
Scherzer O 1995 {\em Journal of Mathematical Analysis and Applications\/} {\bf 194} 911--933

\bibitem{yao2007early}
Yao Y, Rosasco L and Caponnetto A 2007 {\em Constructive Approximation\/} {\bf 26} 289--315

\bibitem{kukavcka2017regularization}
Kuka{\v{c}}ka J, Golkov V and Cremers D 2017 {\em arXiv preprint arXiv:1710.10686\/}

\bibitem{chen2018neural}
Chen R~T, Rubanova Y, Bettencourt J and Duvenaud D~K 2018 {\em Advances in neural information processing systems\/} {\bf 31}

\bibitem{defrise2011algorithm}
Defrise M, Vanhove C and Liu X 2011 {\em Inverse Problems\/} {\bf 27} 065002

\bibitem{vogel2002computational}
Vogel C~R 2002 {\em Computational methods for inverse problems\/} (SIAM)

\bibitem{german2011nonlinear}
German-Sallo Z 2011 {\em UbiCC J\/} {\bf 6} 895--900

\end{thebibliography}

\providecommand{\newblock}{}

\end{document}